\documentclass[12pt]{amsart}
\usepackage{amssymb} 

\usepackage[cmtip,matrix,arrow]{xy}
\CompileMatrices
\SilentMatrices
\SelectTips{cm}{}
\newdir{((}{{\lhook\kern-.1em}}
\newdir{))}{{\rhook\kern-.1em}}
\newdir{ >}{{}*!/-5pt/@{>}}

\newcommand{\A}{\mathcal{A}}
\newcommand{\bt}{\overline{\theta}}

\newcommand{\Co}{\mathcal{C}}
\newcommand{\D}{\mathbf{D}}
\newcommand{\n}{\mathbf{n}}

\newcommand{\bfA}{\mathbf{A}}
\newcommand{\Lg}{L\alg{g}}
\newcommand{\p}{\partial}
\newcommand{\OG}{\Omega G}

\newcommand{\Gau}{\mathcal{G}}

\newcommand{\C}{\mathbb{C}}
\newcommand{\R}{\mathbb{R}}

\newcommand{\alg}[1]{\mathfrak{#1}}

\newcommand{\ve}{\varpi}
\newcommand{\gen}{\chi}
\newcommand{\curv}{\on{curv}}

\DeclareMathOperator{\Ad}{Ad}

\DeclareMathOperator{\Hol}{Hol}

\DeclareMathOperator{\im}{Im}

\DeclareMathOperator{\diag}{diag}
\newcommand\lie[1]{\mathfrak{#1}}

\newcommand{\h}{\lie{h}}
\newcommand{\g}{\lie{g}}
\newcommand{\z}{\lie{z}}
\renewcommand{\a}{\lie{a}}
\renewcommand{\t}{\lie{t}}

\newcommand\Sig{\Sigma}

\newcommand\sig{\sigma}

\newcommand\Om{\Omega}
\newcommand\om{\omega}
\newcommand{\Del}{\Delta}

\newcommand{\f}{\frac}
\newcommand{\Z}{\mathbb{Z}}
\newcommand{\fus}{\circledast} 
\newcommand{\Alc}{\lie{A}} 
\renewcommand{\l}{\langle}
\renewcommand{\r}{\rangle}

\newcommand{\hra}{\hookrightarrow}

\newcommand{\ol}{\overline}
\newcommand{\olt}{\ol{\theta}}
\newcommand{\bib}{\bibitem}
\newcommand{\on}{\operatorname}
\newcommand{\Diff}{\on{Diff}}
\renewcommand{\d}{{\mbox{d}}}
\newcommand{\ti}{\tilde}
\newcommand\qu{/\kern-.7ex/} 

\newcommand{\M}{\mathcal{M}}

\theoremstyle{plain}
\newtheorem{theorem}{Theorem}[section]
\newtheorem{lemma}[theorem]{Lemma}
\newtheorem{corollary}[theorem]{Corollary}
\newtheorem{prop}[theorem]{Proposition}

\theoremstyle{definition}
\newtheorem{definition}[theorem]{Definition}

\newtheorem{remark}{Remark}[section]

\newtheorem{example}{Example}[section]

\begin{document}
\title[Lie group valued moment maps]{Lie group valued moment maps}

\author{Anton Alekseev}
\address{Institute for Theoretical Physics \\ Uppsala University \\
Box 803 \\ \mbox{S-75108} Uppsala \\ Sweden}
\email{alekseev@teorfys.uu.se}

\author{Anton Malkin}
\address{Department of Mathematics \\ Yale University \\ New Haven \\
CT 06520 \\ USA}
\email{malkin@math.yale.edu}

\author{Eckhard Meinrenken}
\address{Massachusetts Institute of Technology, Department of Mathematics,
Cambridge, Massachusetts 02139 \\ USA }
\email{mein@math.mit.edu}

\date{\today}

\begin{abstract}
We develop a theory of ``quasi''-Hamiltonian $G$-spaces for which the
moment map takes values in the group $G$ itself rather than in the
dual of the Lie algebra. The theory includes counterparts of
Hamiltonian reductions, the Guillemin-Sternberg symplectic
cross-section theorem and of convexity properties of the moment map.
As an application we obtain moduli spaces of flat connections on an
oriented compact 2-manifold with boundary as quasi-Hamiltonian
quotients of the space $G^2\times \cdots\times G^2$.
\end{abstract}


\maketitle


\sloppy

\section{Introduction}

The purpose of this paper is to study Hamiltonian group actions 
for which the moment map takes values not in the dual of the Lie algebra 
but in the group itself. 

For the circle group $S^1$ this situation 
has been studied in literature (see e.g. \cite{McD,Wm}). 
The standard example is the real 2-torus $T^2=S^1\times S^1$, with its
standard area form and the circle acting by rotation of the
first $S^1$; the moment map is given by projection to the second $S^1$.
It is in fact known 
that every symplectic $S^1$-action on a symplectic manifold 
for which the 2-form has integral cohomology class admits an 
$S^1$-valued moment map.

In this paper we consider group valued moment maps for general
non-abelian 
compact Lie groups $G$. One theory encountering group valued
moment maps is the theory of Poisson-Lie group actions on symplectic
manifolds \cite{L}, \cite{LW}, where $G$ is a Poisson Lie group and 
the target of the moment map is
the dual Poisson-Lie group $G^*$. In \cite{A} it was shown that for
compact, connected, simply connected Lie groups $G$ this theory is
equivalent to the standard theory of Hamiltonian actions.

In this paper we introduce the notion of ``quasi''-Hamiltonian
(q-Hamiltonian) $G$-spaces consisting of a $G$-manifold $M$, an
invariant 2-form $\om$ and a group valued moment map $\mu:\,M\to G$
satisfying certain natural compatibility conditions. It turns out that
for non-abelian $G$ these spaces differ in many respects from
Hamiltonian $G$-spaces. In particular, the conditions that the 2-form
$\om$ be non-degenerate and closed have to be replaced by somewhat
more complicated conditions.  In spite of these differences,
Hamiltonian reductions of these spaces are defined and result in
spaces with symplectic forms.

Basic examples for q-Hamiltonian $G$-spaces are conjugacy
classes in $G$. Another example is the ``double'' $D(G)=G\times G$
generalizing the above $T^2$-example. Note that for
$G$ compact and simply-connected, 
$D(G)$ does not admit a symplectic structure because its 
second cohomology is trivial.
Yet, it admits a (minimally degenerate) q-Hamiltonian structure. 

As an application we obtain the moduli space $M(\Sig)$ of flat connections
on a closed 2-manifold $\Sig$ of genus $k$ as a q-Hamiltonian quotient
of the space $G^{2k}$. Our construction is a reinterpretation
of the construction 
due to Jeffrey\cite{J1} and Huebschmann\cite{H} (see also \cite{GHJW})
which represents 
$M(\Sig)$ as a symplectic quotient of a certain finite dimensional 
non-compact symplectic
space $X$ (which is in fact  an open subset of $G^{2k}$). 
Their   construction
is based on the group cohomology approach of Goldman\cite{G},
Karshon\cite{K} and Weinstein\cite{W}. For another
 closely related construction see King-Sengupta\cite{KS}.

We show that the space $G^{2k}$ admits a q-Hamiltonian structure
with a $G$-valued moment map corresponding to the $G$-action by
simultaneous conjugations. Then $M(\Sig)$ is obtained as a
q-Hamiltonian quotient of $G^{2k}$.  The closedness and non-degeneracy
of the 2-form on the moduli space follows from the basic properties of
the reduction procedure.  More generally, if $\Sig$ has a boundary and
$\Co=\{\Co_j\}$ are conjugacy classes of holonomies associated to the
boundary components, the corresponding moduli space $M(\Sig,\Co)$ is a
q-Hamiltonian reduction of the space $G^{2(r+k)}$ where $r+1$ is the
number of boundary components. We will explicitly describe the 2-form
on $G^{2(r+k)}$ that gives rise to the symplectic form on moduli space
and check that the answer coincides with Atiyah-Bott's \cite{AB} gauge
theoretic construction of the symplectic form on $M(\Sig,\Co)$.

The paper is organized as follows. The definition of a
q-Hamiltonian $G$-space with $G$-valued moment map is given in
Section 2. Section 3 contains basic examples of q-Hamiltonian
$G$-spaces and Section 4 discusses some of their basic properties. In
Section 5 we show that Hamiltonian reduction extends to the present
setting. In Section 6 we define the ``fusion product'' of two
$G$-Hamiltonian spaces, with moment map the pointwise product of the
two moment maps. In Section 7 we prove a q-Hamiltonian version of the
Guillemin-Sternberg symplectic cross-section theorem and discuss
convexity properties of the moment map. In Section 8 we explain the
relation of our theory to Hamiltonian loop group spaces, we prove that
there is a natural one-to-one correspondence between compact
q-Hamiltonian $G$-spaces and Hamiltonian $LG$-spaces with proper
moment map. Section 9 contains the application of our results to
moduli spaces of flat connections. In Section 10 we explain the
relation to the Lu-Weinstein theory of Poisson-Lie group actions on
symplectic manifolds.\\[.2in]

 \noindent{\bf Acknowledgements:}\\ A.A. is grateful to Department of
Physics, Yale University and to Institute for Theoretical Physics, ETH
Z\"{u}rich for hospitality during the period when this paper was
written.  During his stay at Yale the research of A.A. was supported
in part by DOE Outstanding Young Investigator Award. A.A. thanks
B. Enriquez, G. Felder, J.-C. Hausmann and A. Knutson for useful
discussions.

E.M. was supported by a Feodor Lynen Fellowship from the Humboldt 
Foundation. He is grateful to C. Woodward for 
many interesting suggestions.

\section{Quasi-Hamiltonian $G$-spaces}\label{QHG}
In this section we recall the definition of Hamiltonian $G$-spaces
with $\g^*$-valued moment maps and then present our definition for 
``quasi''-Hamiltonian $G$-spaces with $G$-valued moment maps. 

\subsection{Hamiltonian $G$-spaces}\label{HG}
Throughout this paper $G$ denotes a compact Lie group with Lie algebra
$\alg{g}$. A $G$-manifold is a manifold $M$ together with an action
$\A: G\times M\rightarrow M$.  Given $g\in G,\,x\in M$ we will often
write
\begin{equation*}
\A(g, x)=x^g,
\end{equation*}
and for $\xi\in \alg{g}$ we denote by $v_\xi$ the generating vector field 
on $M$. Given a  closed invariant 2-form $\om$ on a $G$-space 
$(M,\A)$, the contraction $\iota(v_\xi)\om$ is closed because 
$$\d\iota(v_\xi)\om=\mathcal{L}_{v_\xi}\om-\iota(v_\xi)\d\om=0.$$
In symplectic geometry, one is mainly interested in the case that 
$\iota(v_\xi)\om$ is exact and $\om$ is non-degenerate:
\begin{definition}
A Hamiltonian $G$-space $(M,\A,\om,\mu)$
is a $G$-manifold $(M,\A)$, together with 
an invariant 2-form $\omega\in\Om^2(M)^G$ 
and an equivariant {\em moment map} $\mu\in C^\infty(M,\alg{g}^*)^G$  
such that:
\begin{itemize}
\item[(A1)]\label{A1}
The form $\omega$ is closed: $\d\om=0$. 
\item[(A2)]\label{A2}
The moment map satisfies
\begin{equation*}  
\iota(v_{\xi}) \omega = \d \l \mu, \xi \r \ \ \mbox{ for all }\ \
\xi\in \alg{g}.
\end{equation*}
\item[(A3)]\label{A3}
The form $\omega$ is non-degenerate.
\end{itemize}
\end{definition} 
Simple consequences of the axioms are the following description 
of the kernel and image of the derivative of the 
moment map: 
\begin{equation}\label{eq:kernelimage}
\on{im}(\d_x\mu)=\g_x^0,\ \ \ker(\d_x\mu)^\om=\{v_\xi(x),\xi\in\g\}.
\end{equation}
Here $\g_x$ is the isotropy algbra of $x$ and $\g_x^0$ its annihilator 
in $\g^*$, and for any subspace $E\subset T_xM$ the subspace 
$E^\om:=\{v\in T_xM,\,\om(v,w)=0\mbox{ for all }w\in E\}$
is its $\om$-orthogonal complement. 

\subsection{q-Hamiltonian $G$-spaces}\label{sec:definiton}

Let us try to develop the notion of a ``quasi''-Hamiltonian $G$-space
$(M,\A,\om,\mu)$ with a $G$-valued moment map $\mu:\,M\to G$.  We
denote by $(\cdot,\cdot)$ some choice of an invariant positive 
definite inner product on
$\g$ which we use to identify $\g\cong \g^*$, and by
$\theta,\olt\in\Om^1(G,\g)$ the left- and right- invariant
Maurer-Cartan forms. (In a faithful matrix representations for $G$,
$\theta=g^{-1}\d g$ and $\olt=\d g\,g^{-1}$.) If $G$ is abelian (so that 
$\theta=\olt$), the natural 
replacement for the moment map condition is 
\begin{equation}\label{firsttry}
\iota(v_\xi)\om=\mu^*(\theta,\xi).
\end{equation}
If $G$ is non-abelian, condition \eqref{firsttry} does not work 
as 
it is incompatible with the anti-symmetry of $\om$. One 
is forced to replace it by 
\begin{equation}\label{secondtry}
\iota(v_\xi)\om=\f{1}{2}\mu^*(\theta+\olt,\xi).
\end{equation}
If we want $\om$ to be $G$-invariant 
this condition is no longer compatible with $\d\om=0$: 
Indeed, 
\begin{equation}\label{closedcond}
0\stackrel{!}{=}\mathcal{L}_{v_\xi}\om=(\d \iota(v_{\xi}) + \iota(v_{\xi}) \d) \omega=
\frac{1}{2}\mu^*\d \,(\theta+\bt, \xi) +\iota(v_{\xi}) \d \omega. 
\end{equation}
This equation can be rewritten as follows. Let 
$\gen\in\Om^3(G)$ denote the canonical closed bi-invariant 3-form on $G$:
\begin{equation}\label{generator}
\gen = \frac{1}{12} (\theta, [\theta, \theta]) 
=\frac{1}{12} (\bt, [\bt, \bt]).
\end{equation}
Denote by $v_\xi^r$ and $v_\xi^l$ the right and left invariant vector 
field on $G$ generated by $\xi$. The fundamental vector field for the adjoint 
action is $v_\xi=v_\xi^r-v_\xi^l$
so that 
\begin{equation}\label{Fundvf}
\iota(v_\xi)\theta=\Ad_{g^{-1}}\,\xi-\xi,\,\,\iota(v_\xi)\olt=\xi-\Ad_g\xi
\end{equation}
which together with the structure equations
$\d\theta=-\f{1}{2}[\theta,\theta]$ and $\d\olt=\f{1}{2}[\olt,\olt]$
gives
\begin{equation}\label{contraction}
\iota(v_\xi) \gen
=\f{1}{2}\,\d(\olt+\theta,\xi).
\end{equation}
Using these formulas condition \eqref{closedcond} becomes 
$$ 0\stackrel{!}{=}\iota(v_\xi)(\d\om+\mu^*\gen)$$
so that we are lead to require $\d\om=-\mu^*\gen$. 
However, the moment map condition \eqref{secondtry} is 
in general also incompatible with non-degeneracy of $\om$ since it implies that  
all generating vectors $v_\xi(x)$ with $\xi$ a solution of 
$\Ad_{\mu(x)}\xi=-\xi$ have to lie in the kernel of $\om_x$. 

We are therefore
lead to the following ``minimal'' definition.
\begin{definition}\label{def:qham}
A quasi-Hamiltonian $G$-space is a $G$-manifold $(M,\A)$ together with 
an invariant 2-form $\omega\in\Om(M)^G$ and 
an equivariant map
$\mu\in C^\infty(M,G)^G$ such that: 
\begin{enumerate}
\item[(B1)]\label{B1}
The differential  of $\omega$ is given by:
\begin{equation*}  
\d\omega=  - \mu^*\gen.
\end{equation*}
\item[(B2)]\label{B2}
The map $\mu$ satisfies
\begin{equation*} 
\iota(v_{\xi}) \omega =\f{1}{2} \mu^* (\theta+\bt , \xi).
\end{equation*}
\item[(B3)]\label{B3}
At each $x\in M$, the kernel of $\omega_x$ is given by 
\begin{equation*} 
\ker \omega_x= \{ v_{\xi},\,\xi\in\ker(\Ad_{\mu(x)}+1) \} .
\end{equation*}
\end{enumerate}
We will refer to $\mu$ as a moment map.       
\end{definition}

In the following section we will give some examples of 
q-Hamiltonian $G$-spaces. Let us, however, first make a remark
on the definition. 

\begin{remark}\label{rem:equivcoh}
For ordinary Hamiltonian $G$-spaces, the defining conditions 
can be expressed elegantly in terms of the de Rham model for 
equivariant cohomology. Let $\Om^*_G(M)=\bigoplus_k \Om^k_G(M)$ 
be the complex of equivariant differential forms,
\begin{equation*}
\Om^k_G(M)=\bigoplus_{2l+j=k}(\Om^j(M)\otimes S^l\g^*)^G
\end{equation*}
with differential 
\begin{equation*}
(\d_G\alpha)(\xi)=\d(\alpha(\xi))-\iota(v_\xi)\alpha(\xi).
\end{equation*}
Conditions (A1) and (A2)
can be summarized by the requirement that $\om_G(\xi)=\om+\l \mu, \xi \r $ 
be an equivariantly closed form in $\Om^2_G(M)$:
\begin{equation} \label{tdto}
\d_G \om_G(\xi)= \d \om + \d \l \mu, \xi \r - \iota(v_\xi) \om=0.
\end{equation}
However, one can also take a slightly different point of view
and rewrite \eqref{tdto} as follows
\begin{equation*}
\d_G  \om= \d \om -  \iota(v_\xi) \om = - \d \l \mu, \xi \r = 
- \mu^* \gen_G(\xi).
\end{equation*}
Here $\gen_G\in \Om^3_G(\g^*)$ is an equivariant 3-form defined
as  $\gen_G(\xi)=\d\l a,\xi\r$, where $a:\g^*\to\g^*$ is the identity 
map. 
This latter formulation extends to q-Hamiltonian $G$-spaces;
the relevant $\d_G$-closed equivariant 3-form 
$\gen_G\in\Om^3_G(G)$ is 
\begin{equation*}
\gen_G(\xi)=\gen+\f{1}{2}(\theta+\olt,\xi).
\end{equation*}
\end{remark}

\section{Examples of q-Hamiltonian $G$-spaces} \label{example}
\subsection{Conjugacy classes in $G$} \label{conj}
Basic examples of Hamiltonian $G$-spaces are provided by coadjoint
orbits $\mathcal{O}\subset\g^*$. Their q-Hamiltonian 
counterparts are conjugacy classes $\Co\subset G$.

\begin{prop}[Conjugacy classes]\label{pr:Conj}
For every conjugacy class $\Co\subset G$ 
there exists a unique invariant 2-form $\om\in\Om^2(\Co)^G$ for which 
$\Co$ becomes a q-Hamiltonian $G$-space, with 
moment map the embedding $\mu:\,\Co \hra G$. 
The value of $\om$ at $f\in\Co$ is given 
on fundamental vector fields $v_\xi,\,v_\eta$ by 
\begin{equation}\label{qKKS}
\om_f (v_\xi,v_\eta)=
\f{1}{2}\,\Big((\eta,\Ad_f \xi)-(\xi,\Ad_f \eta)\Big).
\end{equation}
\end{prop}
The 2-form \eqref{qKKS} on conjugacy classes plays an 
important role in \cite{GHJW}.
\begin{proof}
Clearly $\om$ is $G$-invariant. The fact that \eqref{qKKS} satisfies 
condition (B2) is tautological because 
\begin{equation*}
\om_f (v_\xi,v_\eta)=
\f{1}{2}\big(\Ad_{f^{-1}}\eta-\Ad_f\eta,\xi)
=\f{1}{2}\,\iota(v_\eta)\big(\theta_f+\olt_f,\xi).
\end{equation*}
To check 
(B1) 
consider the projection $\pi_f: G\rightarrow \Co$ defined by
$\pi_f(u) = \Ad_u f$.
The pull-back of $\om$ is the left-invariant 2-form
\begin{equation*}
\pi_f^*\om=\f{1}{2}(\Ad_f \theta,\theta).
\end{equation*}
We have 
\begin{equation*}
\pi_f^*\d\om=\d\pi_f^*\om=-\f{1}{2}(\Ad_f [\theta,\theta],\theta)
+\f{1}{2}(\Ad_f \theta,[\theta,\theta])
\end{equation*}
Using 
$\pi_f^*\theta=\Ad_{uf^{-1}}\theta-\olt$
one verifies that this last expression is equal to 
$-\pi_f^*\gen$ which shows (B1). 
Suppose next that $\xi\in \g$ is 
such that $v_\xi$ is in the kernel of $\om_f$. By definition of $\om$ this 
means 
\begin{equation*}
\Ad_f \xi-\Ad_{f^{-1}}\xi=0,
\end{equation*}
or $\xi\in\ker(\Ad_{f^2}-1)$.
The kernel of $(\Ad_{f^2}-1)$ is a direct sum 
$$ \ker(\Ad_{f^2}-1)=\ker(\Ad_{f}-1)\oplus \ker(\Ad_{f}+1).$$
For any $\xi\in \ker(\Ad_{f}-1)$ we have $v_\xi(f)=0$.  This shows
that the kernel of $\om_f$ is given by (B3).  Uniqueness of $\om$ is a
consequence of (B2) since $\Co$ is a homogeneous space.
\end{proof}

\begin{remark}
Since the kernel of the 2-form $\om$ on $\Co$ has constant rank, it
defines a distribution on $\Co$. It is integrable, with leaves given
by the orbits $\Ad(Z_{f^2})\cdot f\subset \Ad(G)\cdot f$ as $f$ ranges over
$\Co$.  In particular, if $f^2$ is contained in the center $Z(G)$ then
the 2-form $\om$ on $\Co$ is identically zero.
\end{remark}

\subsection{Double $D(G)$} \label{double}
Our next example is a $q$-Hamiltonian $G\times G$-space 
which plays the same role as the cotangent
bundle $T^*G$ in the category of Hamiltonian $G$-spaces. 
Since this new
q-Hamiltonian space is a product of two copies of $G$ we often
call it a double of $G$ and denote by $D(G)$ (alluding to
the definition of the double in the theory of Quantum Groups).
Let 
\begin{equation*}
D(G):=G\times G
\end{equation*}
and let $a$ and $b$ be projections to the first and
second factor in the direct product. Introduce a $G\times G$-action
$\A_D$ on $D(G)$ via
\begin{equation*} 
(a,b)^{(g_1,g_2)}=(g_1a g_2^{-1},g_2 b g_1^{-1}).
\end{equation*}
Define a moment map $\mu_D=(\mu_1,\mu_2): D(G)\rightarrow G^2$ as
\begin{equation*}
\mu_1(a,b)= ab,\ \ \mu_2(a,b)=a^{-1}b^{-1}
\end{equation*}
and let the 2-form $\omega_D$ be defined by
\begin{equation*}
\omega_D= \f{1}{2}(a^* \theta,\,b^* \olt) + \f{1}{2}
(a^* \bt,\,b^* \theta) .
\end{equation*}
\begin{prop}[The double $D(G)$]
The quadruple $(D(G), \A_D, \mu_D, \omega_D)$ is a 
q-Hamiltonian $G\times G$-space. 
\end{prop}
\begin{proof}
The equivariance of the moment map is immediate from the definition.
Using
\begin{equation}\label{mupullback}
\mu_1^*\theta=\Ad_{b^{-1}} a^*\theta+b^*\theta,\,\,
\mu_2^*\theta=-\Ad_{b}a^*\olt-b^*\olt
\end{equation}
we compute
\begin{equation*}
\mu_1^* \gen = a^* \gen +b^* \gen -\f{1}{2}\, d (a^* \theta, b^* \olt),\ \ 
\mu_2^* \gen = -a^* \gen - b^* \gen -\f{1}{2}\, d (a^* \olt, b^* \theta)
\end{equation*}
which gives the required property $d\om_D=-\mu^*\gen$. 
Next, let $\xi=(\xi_1,\xi_2)\in\g \oplus\g$ 
and $v_\xi$ the corresponding vector field on $D(G)$. 
If $\xi_2=0$ we find, using \eqref{mupullback}, 
\begin{equation*}
\iota(v_\xi)\om_D=\f{1}{2}(\,a^*\olt+b^* \theta +
\Ad_{b^{-1}}a^*\theta+\Ad_{a}b^* \olt ,\xi_1)=
\f{1}{2}\mu_1^*(\theta+\olt,\xi_1).
\end{equation*}
A similar calculation gives the condition for $\xi_1=0$. 

Suppose now that $v\in T_{(a,b)}D(G)$ is in the kernel of $\om$. 
Let $(\xi,\eta)$ be such that $v$ is the value at $(a,b)$ of 
the right-invariant vector field
$v=(v_\xi^r,v_\eta^r)$. The condition 
\begin{equation*}
0=\iota(v)\om_D|_{(a,b)}=
(\Ad_a b^*\olt+b^* \theta,\,\xi)+(\Ad_{b^{-1}} a^*\theta+a^*\olt,\,\eta)
\end{equation*}
gives two equations 
\begin{equation*}
\xi+  \Ad_{ab}   \xi=0,\,\ 
\eta+\Ad_{ab} \eta=0.
\end{equation*}
Setting  
\begin{equation*}
\alpha=\f{1}{2}(\xi+\eta)\,\ \ \beta=-\f{1}{2}\Ad_{b}(\xi-\eta).
\end{equation*}
we obtain
\begin{equation}\label{fund}
\xi = \alpha + \Ad_a \beta \ , \ \eta=\alpha+ \Ad_{b^{-1}} \beta. 
\end{equation}
Equation \eqref{fund} says that $v$ is the value at $(a,b)$ of the fundamental vector 
field for $(\alpha,\beta)$. Moreover, 
\begin{equation*}
\alpha+ \Ad_{ab} \alpha=0 \ , \ \beta + \Ad_{a^{-1}b^{-1}} \beta=0.
\end{equation*}
which shows that (B3) is satisfied.
\end{proof}

\begin{remark}\label{rem:othercoord}
We will also make use of different coordinates on $D(G)$,
corresponding to the left trivialization of $T^*G$: setting 
$u=a$ and $v=ba$ the action reads 
\begin{equation*}
(u,v)^{g_1,g_2}=(g_1\,u\,g_2^{-1},\ \Ad_{g_2} v),
\end{equation*}
the moment map is 
\begin{equation*}
\mu_1(u,v)=\Ad_u v,\ \mu_2(u,v)=v^{-1}
\end{equation*}
and the 2-form is given by
\begin{equation*} 
\omega_D= \frac{1}{2} (\Ad_v u^*\theta, u^*\theta)+
\frac{1}{2} (u^*\theta, v^*\theta+v^*\bt ).
\end{equation*} 
\end{remark}

\subsection{From Hamiltonian $G$-spaces 
to q-Hamiltonian $G$-spaces}\label{SHamqHam}
In this section we show how to construct a q-Hamiltonian $G$-space
from a usual Hamiltonian $G$-space 
$(M,\A,\sigma,\Phi)$.
The basic Lemma is (see also Jeffrey \cite{J1})
\begin{lemma}\label{lem:ve}
For $s\in\R$ let $\exp_s:\,\g\to G$ be defined by
$\exp_s(\eta)=\exp(s\eta)$.
The 2-form on the Lie algebra $\g$ given by 
\begin{equation*}
\ve=\f{1}{2}\int_0^1 (\exp_s^*\olt,\f{\p}{\p s}\exp_s^*\olt) \d s
\end{equation*}
is $G$-invariant and 
satisfies $\d\ve=-\exp^*\gen$. If $v_\xi$ is a fundamental vector 
field for the adjoint $G$-action on $\g$ we have 
\begin{equation}\label{iotavixi}
\iota(v_\xi)\ve= -\d (\cdot,\xi)+\f{1}{2}\exp^*(\theta+\olt,\xi).
\end{equation}
\end{lemma}
We omit the proof at this point since Lemma \ref{lem:ve} is 
a consequence of Proposition \ref{veprop} below.

\begin{prop}\label{HamqHam}
Let $(M,\A,\sigma,\Phi)$ be a Hamiltonian $G$-space. Then $M$ 
with 2-form
$\om=\sigma+\Phi^*\ve$ and moment map $\mu=\exp(\Phi)$ satisfies all 
axioms of a q-Hamiltonian $G$-space except possibly 
the non-degeneracy condition, (B3). If the differential $\d_\xi\exp$ 
is bijective for all $\xi\in \Phi(M)$,  
(B3) is satisfied as well and $(M,\A,\om,\mu)$ is 
a q-Hamiltonian $G$-space.
\end{prop}

\begin{proof}
Equivariance of $\mu$ is clear and invariance of $\om$ follow 
from equivariance of $\exp$ and $\Phi$ and invariance of $\ve$.
Condition (B1) follows from  
$$\d\om=\d\sig+\d\Phi^*\ve=-\Phi^* \exp^* \gen=-\mu^*\gen,$$
and (B2) from the calculation 
\begin{equation*}
\iota(v_\xi)\om=\d(\Phi,\xi)+\f{1}{2}\Phi^*\exp^*(\theta+\olt,\xi)
-\d(\Phi,\xi)=\f{1}{2}\mu^*(\theta+\olt,\xi).
\end{equation*}
To check the non-degeneracy condition (B3) 
suppose that $v\in T_xM$ is in the kernel of $\om$. Then 
\begin{equation}\label{vsigma}
\iota(v)\sigma_x=-\iota(v)(\Phi^*\ve)_x.
\end{equation}
Since the 1-form on the right hand side annihilates the kernel 
of $d_x\Phi$, this equation implies that $v$ is $\sigma$-orthogonal 
to $\ker(d_x\Phi)$. By \eqref{eq:kernelimage} the 
$\sig$-orthogonal complement to $\ker(d_x\Phi)$ span of 
fundamental vector fields $v_\xi(x)$, $\xi\in\g$.
Letting $v=v_\xi(x)$ and using (B2) we arrive at the condition
\begin{equation}\label{kernelcondition}
\Phi^*\exp^*(\theta+\olt,\xi)_x=0
\end{equation}
at $x$.
Pairing with a fundamental vector field $v_\eta(x)$, 
we find 
\begin{equation*}
\big(\eta,\big(\Ad_{\mu(x)^2}-1)\xi\big)=0
\end{equation*}
for all $\eta$, which shows 
$$ \xi\in \ker(\Ad_{\mu(x)^2}-1)=\ker(\Ad_{\mu(x)}-1)\oplus
\ker(\Ad_{\mu(x)}+1).
$$
As we remarked 
in section \ref{sec:definiton}, solutions 
of $\Ad_{\mu(x)}\xi=-\xi$ lead to elements in the kernel. 
Therefore it suffices to consider the case 
$\Ad_{\mu(x)}\xi=\xi$, where the condition 
\eqref{kernelcondition} reads
\begin{equation*}
(\Phi^*(\exp^*\theta,\xi))_x=0.
\end{equation*}
If $\Phi(x)\in\g^*\cong\g$ is not a singular value of $\exp$ this
equation says that $\xi$ annihilates the image of the tangent map
$\d_x\Phi$. By \eqref{eq:kernelimage} this means $\xi\in\g_x$, that is
$v_\xi(x)=0$, and the proof is complete.
\end{proof}


\begin{remark}\label{rem:qhamham}
Suppose conversely that $(M,\A,\om,\mu)$ is a q-Hamiltonian 
$G$-space. Assume also that there exists $U\subset \g$ 
such that $\exp$ is 
a diffeomorphism from $U$ onto some subset $V\subset G$ 
containing $\mu(M)$, and let $\log:\,V\to U$ be the inverse.
Reversing the argument in the proof of Proposition 
\ref{HamqHam}
we see that $(M,\A,\om-\mu^* \log^* \ve,\log(\mu))$ is a Hamiltonian 
$G$-space in the usual sense.
\end{remark}

\section{Properties of q-Hamiltonian $G$-spaces}
The following result summarizes a number of consequences of 
Definition \ref{def:qham}, analogous to \eqref{eq:kernelimage}
for Hamiltonian $G$-spaces. 

\begin{prop}\label{KernelIsomorphism}
Let $(M,\A,\om,\mu)$ be a q-Hamiltonian $G$-space and
$x\in M$. For any subspace $E\subset T_xM$ let $E^\om=\{v\in T_xM,\,\om(v,w)=0
\mbox{ for all }w\in E\}$ denote its $\om$-orthogonal complement. 
\begin{enumerate}
\item
The map \ \  
$\ker(\Ad_{\mu(x)}+1)\to \ker\om_x,\,\,\xi\mapsto v_\xi(x)$
is an isomorphism. 
\item
$ \ker(\d_x\mu)\cap\ker\om_x=\{0\}.$
\item
$ \on{im}(\mu^*\theta)_x=\g_x^\perp.$
\item
$(\ker{\d_x\mu})^\om=\{v_\xi(x),\xi\in\g\}$.
\end{enumerate}
\end{prop}
\begin{proof}
Observe first that there is an orthogonal splitting of $\g$,
\begin{eqnarray*}
\g&=&\ker(\Ad_{\mu(x)}+1) \oplus \on{im}(\Ad_{\mu(x)}+1)\\
&=&\ker(\Ad_{\mu(x)}+1) \oplus\ker(\Ad_{\mu(x)}-1) \oplus
\on{im}(\Ad_{\mu(x)^2}-1) 
\end{eqnarray*}
and that $\g_{\mu(x)}=\ker(\Ad_{\mu(x)}-1)$. \\
1. By the non-degeneracy condition (B3) the map  
is surjective. 
On the other hand if $v_\xi(x)=0$ then $v_\xi(\mu(x))=0$ 
by 
equivariance of the moment map or equivalently
$\xi\in\ker(\Ad_{\mu(x)}-1)\subset \on{im}(\Ad_{\mu(x)}+1)$.
This shows injectivity.\\
2. Let $v\in \ker(\d_x\mu)\cap\ker \om_x$. 
By 1. we can write $v=v_\eta$ with 
$\eta\in\ker (\Ad_{\mu(x)}+1)$. By equivariance of the moment map, 
$0=\d_x\mu(v_\eta(x))=v_\eta(\mu(x))$, which shows $\eta\in\g_{\mu(x)}=
\ker(\Ad_{\mu(x)}-1)$. Thus $\eta=0$ and consequently $v=0$.\\
3.
Using the defining equation for the moment map,
\begin{equation*}
\om(v_\xi,v)=\frac{1}{2}\iota(v)\mu^*(\theta+\olt,\xi)
=\frac{1}{2}\big((\Ad_{\mu}+1)\iota(v)\mu^*\theta,\xi\big)
\end{equation*}
and Property 1. we have
\begin{eqnarray*} 
\on{im}(\mu^*\theta)_x\cap\on{im}(\Ad_{\mu(x)}+1)
&=&(\Ad_{\mu(x)}+1)\on{im}(\mu^*\theta)_x=\{\xi,v_\xi(x)\in\ker\om_x\}^\perp\\
&=&(\g_x\oplus  \ker(\Ad_{\mu(x)}+1))^\perp=\g_x^\perp\cap \on{im}(\Ad_{\mu(x)}+1).
\end{eqnarray*}
On the other hand, 
equivariance of the moment map together with \eqref{Fundvf}
shows that
$$\on{im}(\mu^*\theta)_x\supset \on{im}(\Ad_{\mu(x)^{-1}}-1)=
\ker(\Ad_{\mu(x)}-1)^\perp\supset\ker(\Ad_{\mu(x)}+1).$$ 
Since
$\g_x\subset\g_{\mu(x)}=\ker(\Ad_{\mu(x)}-1)\subset \on{im}(
\Ad_{\mu(x)}+1)$ 
so that also $\g_x^\perp\supset\ker(\Ad_{\mu(x)}+1)$
these two equations prove 3.\\
4. The inclusion $\supset$ is a direct consequence from the defining property
(B2) of the moment map. Equality follows by dimension count:
Using 2. and 3., we have 
$$\dim(\ker\d_x\mu)^\om=\dim M-\dim\ker(\d_x\mu)=
\dim \on{im}(\d_x\mu)=\dim \g-\dim \g_x.$$
\end{proof}

\begin{remark}\label{rem:refinement}
We shall also need the following refinement of Property 4.: 
Suppose that $G$ is a product $G=G_1\times G_2$, 
and let $\mu=(\mu_1,\mu_2)$ be the components
of the moment map. Then 
$$ (\ker \d_x\mu_1)^\om=\{v_\xi,\xi\in\g_1\}+ \ker\om_x.$$
To see this note that as a direct consequence of Property 3., 
$\on{im}(\d_x\mu_1)=\g_x^\perp\cap \g_1
=(\g_1)_x^\perp\cap\g_1$ which implies the above equation by dimension 
count.
\end{remark}

\begin{prop}\label{Immersion}
Let $(M_j,\A_j,\om_j,\mu_j)$ ($j=1,2$) be q-Hamiltonian $G$-spaces, and 
$F:\,M_1\to M_2$ an equivariant smooth map such that $F^*\om_2=\om_1$ and 
$F^*\mu_2=\mu_1$. Then $F$ is an immersion.
\end{prop}

\begin{proof}
Let $x\in M$. Since $F^*\om_2=\om_1$ we have $\ker \d_x
F\subset\ker(\om_1)_x$. Since $F^*\mu_2=\mu_1$,
Lemma
\ref{KernelIsomorphism} shows that $\d_xF$ restricts to an isomorphism  
$\ker(\om_1)_x\cong  \ker(\om_2)_{F(x)}$. It follows that 
$\ker \d_x F\cap\ker(\om_1)_x=\{0\}$.
\end{proof}
\begin{corollary}
If $(M,\A,\om,\mu)$ is a q-Hamiltonian $G$-space on which $G$-acts
transitively, the moment map $\mu$ is a covering onto a conjugacy 
class $\Co\subset G$.
\end{corollary}
\begin{proof}
Let $x\in M$ and $f=\mu(x)$. 
By equivariance $\mu$ is a submersion from $M$ onto $\Co=\Ad(G)\cdot f$. 
(B2) shows that $\om$ is the pull-back by 
$\mu$ of the 2-form on $\Co$. Therefore $\mu$ is an immersion by 
Proposition \ref{Immersion}.
\end{proof}

\begin{prop}[Inversion]\label{Mminus}
If $(M,\A,\om,\mu)$ is a q-Hamiltonian $G$-space then the quadruple
$(M,\A,-\om,\mu^{-1})$, 
is also a q-Hamiltonian $G$-space. 
\end{prop}

\begin{proof}
This is an immediate consequence of the fact that under the inversion map 
$\on{Inv}:G\to G,\,g\mapsto g^{-1}$
\begin{equation*}
\on{Inv}^*\theta=-\olt,\,\on{Inv}^*\olt=-\theta,\,\on{Inv}^*\gen=-\gen.
\end{equation*}
\end{proof}

We denote the q-Hamiltonian $G$-space given in Proposition \ref{Mminus}
by $M^-$. Note that if $\Co$ is the conjugacy class of $f\in G$ then 
$\Co^-$ is the conjugacy class of $f^{-1}$.

All of the above results are natural analogues of well-known 
facts about Hamiltonian $G$-spaces. The following Theorem
introduces a nontrivial automorphism of q-Hamiltonian $G$-spaces 
that does not have a counterpart for Hamiltonian $G$-spaces.
As we shall see later, this automorphism
corresponds to Dehn twists  of 2-manifolds.

\begin{theorem}[Twist automorphism] \label{Dehn}
Let $(M,\A,\om,\mu)$ be a q-Hamiltonian $G$-space.
The map 
\begin{equation*}
Q:\,M\to M,\,x\mapsto x^{\mu(x)}
\end{equation*}
is a diffeomorphism satisfying $Q^*\om=\om,\,Q^*\mu=\mu$.
\end{theorem}

\begin{proof}
By equivariance of the moment map, 
$Q^*\mu=\Ad(\mu)\mu=\mu$. The map 
$x\mapsto x^{\mu(x)^{-1}}$ is an inverse to $Q$ so that $Q$ is a 
diffeomorphism. The tangent map 
to $Q$ is given by 
\begin{equation*}
\d_x Q(v)=\d_x \A_{\mu(x)}v+v_\xi(Q(x))
\end{equation*}
where $\xi= \iota(v)(\mu^*\olt)_x\in\g$. Letting $v_j\in T_x M$
be two tangent vectors and  $\xi_j=\iota(v_j)(\mu^*\olt)_x$, we 
find that the expression
$Q^*\om_x(v_1,v_2)=\om_{Q(x)}(\d_x Q(v_1),\d_x Q(v_2))$ 
is the sum of the following four terms:
\begin{equation*}
\om_{Q(x)}(\d_x \A_{\mu(x)}v_1,\d_x \A_{\mu(x)}v_2)=\om_x(v_1,v_2)
\end{equation*}
plus
\begin{eqnarray*}
\f{1}{2}\om_{Q(x)}(v_{\xi_1}(Q(x)),\d_x \A_{\mu(x)}v_2)&=&
\f{1}{2}\iota(\d_x \A_{\mu(x)} v_2)\mu^*(\theta+\olt,\xi_1)_{Q(x)}\\
&=&\f{1}{2}\iota(v_2)\mu^*(\theta+\olt,\Ad(g^{-1})\xi_1)_{x} \nonumber \\ 
&=&\f{1}{2}(\Ad(\mu(x)^{-1})\xi_2  -\xi_2,\Ad(\mu(x)^{-1})\xi_1)\nonumber \\
&=&\f{1}{2}(\xi_1,\xi_2 -\Ad(\mu(x))\xi_2)\nonumber 
\end{eqnarray*}
minus the same term with $v_1,v_2$ exchanged, plus 
\begin{equation*}
\f{1}{2}\iota\big(v_{\xi_2}(Q(x))\big)\mu^*(\theta+\olt,\xi_1)_{Q(x)}
=\f{1}{2}(\xi_1,\Ad(\mu(x)^{-1})\xi_2 -\Ad(\mu(x))\xi_2).
\end{equation*}
Adding up all contributions we find
$Q^*\om_x(v_1,v_2)=\om_x(v_1,v_2)$.
\end{proof}

\begin{remark}
If $G$ is a product $G=G_1\times G_2$ one has a twist automorphism for 
every factor: That is, both maps $x\mapsto x^{\mu_j(x)}$ are equivariant 
diffeomorphisms preserving 
$\mu$ and $\om$.
\end{remark}

On q-Hamiltonian $G$-spaces one can define $G$-invariant Hamiltonian
vector fields. 

\begin{prop}[Hamiltonian dynamics] \label{HamDyn}
Let $(M,\A,\om,\mu)$ be a q-Hamiltonian $G$-space. 
For every $G$-invariant function 
$F\in C^\infty(M,\R)^G$ there is a unique smooth vector field $v_F$ satisfying
the following conditions:
\begin{equation} \label{Ham}
\iota(v_F) \omega = dF \ , \ \iota(v_F) \mu^* \theta = 0.
\end{equation}
The vector field $v_F$ is $G$-invariant and preserves $\om$ and $\mu$.
\end{prop}

\begin{proof}
By Proposition \ref{KernelIsomorphism}, Property 1 the map 
\begin{equation*}
A:\,TM\to T^*M\oplus \g,\ v\mapsto (\iota(v)\om,\iota(v)\mu^*\theta). 
\end{equation*}
is injective. Its image defines a smooth sub-bundle 
$E\subset T^*M\oplus\g$. We need to show that 
$x\mapsto (\d_x F,0)$ defines a section 
of $E$, the corresponding vector field $v_F$ is then just the 
pre-image under $A$. 
Since $F$ is $G$-invariant, $\d_x F$ annihilates the 
space $\{v_\xi(x),\xi\in\g\}\subset T_x M$ and in particular 
the kernel $\ker \om_x$. 
This shows that $\d_xF$ is contained 
in the image of the map $T_xM\to T_x^*M,\,v\mapsto \iota(v)\om_x$.
Let $v_0\in T_xM$ with $\d_xF=\iota(v_0)\om_x$. 
and set $\xi:=\iota(v_0)(\mu^*\theta)_x$. For all 
$\eta\in \g$ we have, using once again that $F$ is invariant,  
\begin{equation*}
0=\iota(v_{\eta}) \d_xF= - \iota(v_0) \iota(v_{\eta}) \omega_x=
- \frac{1}{2} \iota(v_0) \mu^*(\theta + \bt, \eta)_x = 
- \frac{1}{2} (\xi + \Ad_{\mu(x)^{-1}} \xi, \eta)
\end{equation*}
which shows $\xi\in \ker(\Ad_{\mu(x)}+1)$. Consequently 
$v_1:=v_{\xi}(x)\in \ker \om_x$ and $v=v_0+\f{1}{2}\,v_1$ still 
solves $\iota(v)\om_x=\d_xF$. Using (B2) we have 
\begin{equation*} 
\iota(v)(\mu^*\theta)_x=
\iota(v_0+\f{1}{2}\,v_1)(\mu^*\theta)_x=\xi+\f{1}{2}(\Ad_{\mu^{-1}}\xi-\xi)
=\xi-\xi=0
\end{equation*}
so that $v$ solves $A(v)=(\d_xF,0)$. This shows that 
$(\d_xF,0)$ is in the image of $A$.

$G$-invariance of $v_F$ follows by the $G$-invariance of its 
defining equations. The equation
\begin{equation*}
\mathcal{L}_{v_F}\om= (d\iota(v_F)+\iota(v_F)d) \omega = d dF- 
 \iota(v_F) \mu^* \gen=0 
\end{equation*} 
shows that
the 2-form $\omega$ is $v_F$-invariant. Invariance of $\mu$ 
is equivalent to invariance of the $G$-valued 1-form $\mu^*\theta$:
\begin{equation*}
\mathcal{L}_{v_F}\mu^*\theta=(d\iota(v_F)+\iota(v_F)d)\mu^*\theta
=-\f{1}{2}\,\iota(v_F)\mu^*[\theta,\theta]=0.
\end{equation*}
\end{proof}

\section{q-Hamiltonian reduction} \label{Reduction}

In this Section we show that the usual Hamiltonian
(Meyer-Marsden-Weinstein) reduction procedure can be carried out for
q-Hamiltonian $G$-spaces. We assume that $G$ is a product $G=G_1\times
G_2$ and we consider reductions with respect to the first factor.
Given $f\in G_1$ let $Z_f\subset G_1$ be its centralizer and $\z_f$
the Lie algebra.  Let $M$ be a q-Hamiltonian $G_1\times G_2$-space,
with moment map $(\mu_1,\mu_2)$.  Suppose that $f\in G_1$ is a regular
value of $\mu_1$, so that $\mu_1^{-1}(f)$ is a smooth submanifold.
Proposition \ref{KernelIsomorphism} shows that for all $x\in \mu_1^{-1}(f)$, 
the isotropy algebra 
$(\g_1)_x$ is trivial, 
i.e. $(G_1)_x\subset Z_f$ is a discrete subgroup. 
It follows that the {\em reduced space}
$M_f=\mu_1^{-1}(f)/Z_f$ is an $G_2$-equivariant orbifold.  Not
surprisingly $M_f$ is a q-Hamiltonian $G_2$-space:
\begin{theorem}[q-Hamiltonian reduction]
Let $M$ be a q-Hamiltonian $G_1\times G_2$-space and let 
$f\in G_1$ be a 
regular value of the moment map $\mu_1: M\rightarrow G_1$.
Then the pull-back of the 2-form $\omega$ to $\mu_1^{-1}(f)$ 
descends to the {\em reduced space}
\begin{equation*}
M_f = \mu_1^{-1}(f)/Z_f 
\end{equation*}
and makes it into a q-Hamiltonian $G_2$-space. In particular, if 
$G_2=\{e\}$ is trivial then $M_f$ is a symplectic orbifold. 
\end{theorem}
\begin{proof}
Let $\iota:\,\mu_1^{-1}(f)\hra M$ denote the embedding
and $\pi:\,\mu_1^{-1}(f)\to M_f$ the projection. 
The form
$\iota^* \omega$ is $Z_f\times G_2$-invariant because $\omega$ is 
$G_1\times G_2$-invariant. 
Moreover if $\xi\in \z_f$,
\begin{equation*}
\iota(v_{\xi}) \iota^* \omega = 
\iota^*\iota(v_{\xi}) \omega = 
\iota^* \mu_1^*(\theta+ \bt, \xi) = 0
\end{equation*}
(because $\mu_1\circ \iota:\,\mu_1^{-1}(f)\to \{ f \}$) implies that
$\iota^* \omega$ is $Z_f$-basic. Let $\om_f\in \Om^2(M_f)^{G_2}$ be the
unique 2-form such that $ \pi^* \om_f = \iota^* \omega$.
The restriction $\iota^*\mu_2$ is $Z_f\times G_2$-invariant and descends to an
equivariant map
$(\mu_2)_f\in C^\infty(M_f,G_2)^{G_2}$ satisfying (B2).  
Letting $\gen_1,\,\gen_2$ be the canonical
3-forms for $G_1,G_2$, we have
\begin{equation*}
\pi^*\d\om_f=\d\iota^*\om=\iota^*\d\om=
-\iota^*(\mu_1^*\gen_1+\mu_2^*\gen_2)=
-\iota^*\mu_2^*\gen_2=-\pi^*(\mu_2)_f^* \gen_2
\end{equation*}
so that $\om_f$ satisfies (B1). 
Finally, we need
to check the non-degeneracy condition (B3) 
for $\om_f$. 
The kernel of $\om_f$
at a point $\pi(x)\in M_f$ is just the projection
of $\d_x\pi\ker (\iota^*\om)_x$.
Using Remark \ref{rem:refinement},  
\begin{eqnarray*}
\ker(\iota^*\om)_x&=&\ker(\d_x\mu_1)\cap \ker(\d_x\mu_1)^\om\\
&=&\ker(\d_x\mu_1)\cap\big(\{v_\xi(x),\xi\in\g_1\}+\ker\om_x\big)\\
&=&T_x(Z_f\cdot x)+\{v_\eta(x)|\,\eta\in\ker(\Ad_{\mu_2(x)}+1)\}.
\end{eqnarray*}
\end{proof}
\begin{remark}
The proof shows that if $M$ satisfies the conditions for a 
q-Hamiltonian $G_1\times G_2$-space except for (B3), the 
reduced space $M_f$ is still well-defined and satisfies 
(B1) and (B2). 
\end{remark}

Let us consider a few examples of q-Hamiltonian reduction. 
\begin{example}
\begin{enumerate}
\item For a conjugacy class $M=\Co\subset G$, the reduced space  
$M_f$ is a point if $f\in \Co$, empty otherwise. 
\item
Let $M=D(G)$ be the double defined in the previous section, 
with moment map $\mu_D=(\mu_1,\mu_2)$. Consider the reduction 
with respect to the second $G$-factor, $\mu_2^{-1}(f)/Z_f$. Since 
$\mu_2^{-1}(f)=\{(a,f^{-1}a^{-1})|\,a\in G\}$, with $Z_f$  acting 
diagonally from the right and the first $G$ acting diagonally 
from the left we find that the reduced space is the conjugacy 
class through the element $f^{-1}$.
\item
It follows from the proof that 
if $M$ satisfies (B1) and (B2) but not necessarily (B3), 
the reduced space $M_f$ is still well-defined and satisfies (B1),(B2).
For example let $(M,\A,\sig,\Phi)$ be a Hamiltonian $G$-space, and 
let $(M,\A,\om,\mu)$ be the q-Hamiltonian $G$-space obtained from it,
with $\mu=\exp(\Phi)$ and $\om=\sig+\Phi^*\ve$. Recall that even if 
$\sig$ is non-degenerate $\om$ can fail to satisfy condition (B3). 
Suppose  $f$ is a regular 
value of both $\mu:M\to G$ and $\exp:\g\to G$. Then all  
pre-images $\mu\in\exp^{-1}(f)$ 
are regular values of $\Phi$, and the q-Hamiltonian reduction 
$M_f$ is symplectic and is a disjoint union 
$$ M_f=\coprod_{\exp\mu=f}M_\mu.$$
This follows from the fact that the pull-back of the extra 
term $\Phi^*\ve$ to $\mu^{-1}(f)$ vanishes.
\end{enumerate}
\end{example}

\begin{remark}\label{shiftingtrick}
The q-Hamiltonian structure on $M_f$ can also be obtained as follows. 
Let $\tau$ be the 2-form on the conjugacy class $\Co=\Ad(G)\cdot f$ and 
let $\iota:\,\mu_1^{-1}(\Co)\to M$ and $\pi:\,\mu_1^{-1}(\Co)\to 
\mu_1^{-1}(\Co)/G$ be the embedding and projection. Then 
$$ \pi^*\om_f=\iota^*(\om-\mu^*\tau). $$
\end{remark}

\begin{remark}[Hamiltonian Dynamics commutes with reduction]
The Hamiltonian dynamics defined on $M$ by a $G_1\times G_2$-invariant
Hamiltonian $F$ descends to the reduced space $M_f$. Indeed, the
corresponding vector field $v_F$ descends to $\mu_1^{-1}(f)$ because
it is $G$-invariant and tangent to $\mu_1^{-1}(f)$.  Moreover, the
restriction of the $G_1\times G_2$-invariant function $F$ to
$\mu_1^{-1}(f)$ is $Z_f\times G_2$-invariant and descends to an
$G_2$-invariant function $F_f$ on $M_f$ satisfying
\begin{equation*}
\pi^* \d F_f= \iota^* dF= \iota^* \iota(v_F) \omega  =
\pi^* \iota(\pi (v_{F_f}))\omega_f
\end{equation*}
and 
\begin{equation*}
\pi^*\iota(v_{F_f})(\mu_2)_f^*\theta=\iota(v_F)\iota^*\mu_2^*\theta=0.
\end{equation*}
It follows that
$G_1\times G_2$-invariant  Hamiltonian dynamics commutes with reduction.
\end{remark}

\section{Fusion product} \label{fusion}

In this section we introduce a ring structure on the category of
q-Hamiltonian $G$-spaces.  We call it a fusion product because it
provides a finite-dimensional ``classical analogue'' to fusion
products of representations of quantum groups at roots of unity, and
to fusion products of positive energy representations of loop groups
(see \cite{MW1} and Section \ref{Loopgroup} below).
\begin{theorem}[Fusion product]\label{FusionProduct}
Let $M$ be a q-Hamiltonian $G\times G\times H$-space, 
with moment map $\mu=(\mu_1,\mu_2,\mu_3)$. Let $G\times H$ act by 
the diagonal embedding $(g,h)\to (g,g,h)$. Then $M$ with 2-form 
\begin{equation} \label{Fomega}
\ti{\om}=\om+ \f{1}{2}(\mu_1^* \theta, \mu_2^* \olt)
\end{equation}
and moment map
\begin{equation} \label{Fmoment}
\ti{\mu} = (\mu_1\cdot \mu_2,\mu_3):\,M\to G\times H
\end{equation}
is a q-Hamiltonian $G\times H$-space. 
\end{theorem}

\begin{proof}
The moment map $\ti{\mu}$ is 
equivariant because the group multiplication
is an equivariant map with respect to the action by
conjugations. 
Property (B1) follows from 
\begin{equation}\label{multiplicationidentity}
(g_1g_2)^*\gen=g_1^*\gen+g_2^*\gen-\f{1}{2}\d\,( g_1^*\theta,\,g_2^*\olt).
\end{equation}
thus
\begin{equation*}
d\omega=  
-\mu_1^* \gen - \mu_2^* \gen -\mu_3^*\gen +
\f{1}{2}\d(\mu_1^*\theta,\mu_2^*\olt)= -\ti{\mu}^* \gen.
\end{equation*}
For $\xi\in\g$ and $\eta\in\h$ let 
$v_\xi^1,\,v_\xi^2$ and $v_\eta$ denote the fundamental 
vector fields for the action of the respective factors of 
$G\times G\times H$. The fundamental vector field for the diagonal 
$G$-action is just the sum $v_\xi=v_\xi^1+v_\xi^2$. 
Clearly 
\begin{equation*} 
\iota(v_\eta)\ti{\om}=\iota(v_\eta)\om=\f{1}{2}
\mu_3^*(\theta+\olt,\eta)
\end{equation*}
which verifies (B2) for the $H$-factor. Moreover,
\begin{eqnarray*}
\lefteqn{\iota(v_{\xi}) \ti\om
=  \iota(v_{\xi}^1) \om +
 \iota(v_{\xi}^2) \om + 
\frac{1}{2} 
(\mu_1^* \iota(v_{\xi}) \theta, \mu_2^* \olt) -\frac{1}{2} 
(\mu_1^* \theta, \mu_2^* \iota(v_{\xi}) \olt)} \\ 
&=&
\frac{1}{2} \mu_1^* (\theta +\bt, \xi) + \frac{1}{2} \mu_2^*(\theta + \bt, \xi)
+
\frac{1}{2}  (\Ad_{\mu_1^{-1}} \xi -\xi , \mu_2^*\olt) 
-\frac{1}{2} (\mu_1^* \theta, \xi - \Ad_{\mu_2}\xi)\\&=&\frac{1}{2}  
(\mu_1\mu_2)^* (\theta +\bt, \xi).
\end{eqnarray*}
Finally, we need to check that $\ti\omega$ satisfies the
non-degeneracy condition (B3).
Suppose the vector $v\in T_xM $ 
is in the kernel of $\ti{\om}_x$ (We will omit the basepoint $x$ 
to simplify notation):
\begin{equation}\label{kernel}
0=\iota(v) \ti\om = \iota(v) \om + \frac{1}{2}  (\iota(v) \mu_1^* \theta, \mu_2^* \olt) - \frac{1}{2}
( \mu_1^* \theta, \iota(v) \mu_2^* \olt) .
\end{equation}
Let $\zeta\in\ker(\Ad_{\mu_1}+1)$, so that $v_\zeta^1\in\ker\om_x$.
Contracting \eqref{kernel} with 
$v_\zeta$ we find 
\begin{equation*}
0=\iota(v_{\zeta}^1)\ (\mu_1^*\theta, \iota(v) \mu_2^* \olt)
=((\Ad_{\mu_1^{-1}}-1)\zeta, \iota(v) \mu_2^* \olt)
=-2(\zeta,\iota(v) \mu_2^* \olt).
\end{equation*}
This shows  
$$\iota(v) \mu_2^* \olt\in\ker(\Ad_{\mu_1}+1)^\perp=\on{im}(\Ad_{\mu_1}+1).$$
A similar 
argument applies to $\iota(v) \mu_1^*\theta $. 
We can therefore choose $\xi_1,\xi_2\in\g$ with 
\begin{equation}\label{KeyEquations1}
\iota(v) \mu_2^* \olt=(\Ad_{\mu_1^{-1}}+1)\xi_1,\, \ \ \
\iota(v) \mu_1^*\theta =-(\Ad_{\mu_2}+1)\xi_2
\end{equation}
which turns \eqref{kernel} into 
\begin{gather*}
\iota(v)\om=\f{1}{2}\mu_1^*(\theta+\olt,\xi_1)+
\f{1}{2}\mu_2^*(\theta+\olt,\xi_2)
=\iota(v_{\xi_1}^1)\om +\iota(v_{\xi_2}^2).
\end{gather*}
Changing $\xi_1,\xi_2$ if necessary it follows that 
$v=v_{\xi_1}^1+v_{\xi_2}^2+v_\eta$ for suitable $\xi_1,\xi_2\in\g$
and $\eta\in \ker(\Ad_{\mu_3}+1)$.
Re-inserting this into equations \eqref{KeyEquations1}
we find
\begin{equation*}
(1-\Ad_{\mu_2})\xi_2=(1+\Ad_{\mu_1^{-1}})\xi_1,\,\ \ \
(1+\Ad_{\mu_2})\xi_2=(1-\Ad_{\mu_1^{-1}})\xi_1.
\end{equation*}
Adding these equations gives $\xi_1=\xi_2$, and then either 
equation shows the non-degeneracy condition
\begin{equation*}
\Ad_{\mu_1\mu_2}\xi=-\xi.
\end{equation*} 
\end{proof}

We call the operation of replacing the $G\times G\times H$-action by
the $G\times H$-action on a manifold $M$ internal fusion, and denote
the resulting space by $M_{12}$ (in particular if there are more
$G$-factors involved).  Given two q-Hamiltonian $G\times H_j$-spaces
$M_j$ we define their fusion product $M_1\fus M_2$ to be the
q-Hamiltonian $G\times H_1\times H_2$-space obtained from the
q-Hamiltonian $G\times H_1\times G\times H_2$-space $M_1\times M_2$ by
fusing the $G$-factors.
\begin{remark}
\begin{enumerate}
\item Let $\{\on{pt}\}$ denote the trivial $G$-space, with moment 
map $\on{pt}\mapsto e$. Then $M\fus \{\on{pt}\}=M=\{\on{pt}\}\fus M$ 
for every q-Hamiltonian $G$-space $M$.
\item
Suppose $M$ is a q-Hamiltonian $G\times G\times H$-space. 
Then $(M^-)_{12}=(M_{21})^-$.
\end{enumerate}
\end{remark}

The fusion operation is associative: Given q-Hamiltonian 
$G\times H_j$-spaces $M_j$ we have 
\begin{equation*}
(M_1\fus M_2)\fus M_3= M_1\fus(M_2\fus M_3)
\end{equation*}
More generally
if $M$ is a q-Hamiltonian $G\times G\times G\times H$-space, 
the two q-Hamiltonian $G\times H$-spaces 
$M_{(12)3}$ 
obtained by first fusing the first two $G$-factors 
and $M_{1(23)}$
obtained by first fusing the last two $G$-factors are 
identical. The new 2-form on $M$ in either case is 
given by 
\begin{equation*}
\om+\f{1}{2}(\mu_1^*\theta,\mu_2^*\olt)
+\f{1}{2}(\mu_2^*\theta,\mu_3^*\olt)+\f{1}{2}(\mu_1^*\theta,
\Ad_{\mu_2}\mu_3^*\olt).
\end{equation*}
We shall now show that the fusion product is also commutative
on {\em isomorphism classes} of $q$-Hamiltonian $G$-spaces. 
Switching the two $G$-factors in Theorem \ref{FusionProduct} 
before fusing we obtain a Hamiltonian $G\times H$-space $M_{21}$ with 
the same action, but moment map $(\mu_2\cdot\mu_1,\mu_3)$
and 2-form
\begin{equation*}
\om+\f{1}{2}(\mu_2^*\theta,\mu_1^*\olt).
\end{equation*}
If $G$ is non-commutative  the identity map $M\to M$ does 
not provide an isomorphism. However, we have the following
result. Let $\A^1,\A^2:\,G\to \on{Diff}(M)$ denote the actions
of the two $G$-factors. 
\begin{theorem}[Commutativity of the Fusion Product]
Under the hypotheses of Theorem \ref{FusionProduct} the map
\begin{equation*}
R:\,M\to M,\,x\mapsto \A^2_{\mu_1(x)}(x)
\end{equation*}
is a $G\times H$-equivariant diffeomorphism satisfying 
\begin{gather}
R^*(\mu_2\mu_1)=\mu_1\mu_2,\,R^*\mu_3=\mu_3,\label{RMoment}\\
R^*\Big(\om+\f{1}{2}(\mu_2^* \theta, \mu_1^* \olt)\Big)=
\om+\f{1}{2}(\mu_1^* \theta, \mu_2^* \olt).\label{R2form}
\end{gather}
Thus $R$ gives an isomorphism $R:\,M_{12}\to M_{21}$ of 
q-Hamiltonian $G\times H$-spaces.
\end{theorem}
\begin{proof}
By equivariance of the moment map $R$ is $G\times H$-equivariant, and  
we have 
\begin{equation*}
R^*\mu_1=\mu_1,\,\,R^*\mu_2=\Ad(\mu_1)\mu_2,\,\,R^*\mu_3=\mu_3.
\end{equation*}
proving \eqref{RMoment}. To prove \eqref{R2form} we note that the 
tangent map
to $R$ is 
\begin{equation*}
\d_x R(v)=\d_x\A^2_{\mu_1(x)}(v)+v^2_{\xi}(R(x))
\end{equation*}
where $v^2_\xi$ is the fundamental vector field of 
$\xi:=\iota(v)(\mu_1^*\olt)_{x}$ with respect to $\A^2$. 
By a calculation similar to that in the proof of Theorem \ref{Dehn}
we find 
\begin{eqnarray*}
\lefteqn{
(R^*\om)(v_1,v_2)=\om(v_1,v_2)+\f{1}{2}
\Big(\Ad(\mu_1^{-1})\xi_1, \xi_2+\Ad(\mu_2)\xi_2\Big)}\\
&& -\f{1}{2}\Big(\Ad(\mu_1^{-1})\xi_2, \xi_1+\Ad(\mu_2)\xi_1\Big)
-\f{1}{2}\big(\xi_1,\Ad(\mu_2)\xi_2-\Ad(\mu_2^{-1})\xi_2\big)
\end{eqnarray*}
with $\xi_j=\iota(v_j)(\mu_1^*\olt)$, thus
\begin{equation}\label{R1}
R^*\om=\om+\f{1}{2}(\mu_1^*\theta,\mu_2^*(\theta+\olt)-\Ad(\mu_2)\mu_1^*\theta).
\end{equation}
Furthermore
\begin{eqnarray*} 
R^*\mu_1^*\olt&=&\mu_1^*\olt\\
R^*\mu_2^*\theta&=&-\mu_1^*\olt+\Ad(\mu_1\mu_2^{-1})\mu_1^*\theta+\Ad(\mu_1)
\mu_2^*\theta
\end{eqnarray*}
so that 
\begin{equation}\label{R2}
R^*(\mu_2^*\theta,\mu_1^*\olt)=
(\Ad(\mu_2^{-1})\mu_1^*\theta+\mu_2^*\theta,\mu_1^*\theta)
=(\mu_1^*\theta,\,\Ad(\mu_2)\mu_1^*\theta-\mu_2^*\theta)
.
\end{equation}
Adding \eqref{R1} $+\f{1}{2}$ \eqref{R2} proves the theorem.
\end{proof}
We leave it to the reader to check that the map 
\begin{equation*} 
R':\,M\to M,\,x\mapsto \A^1_{\mu_2(x)^{-1}}(x) 
\end{equation*}
has just the same properties \eqref{RMoment}, \eqref{R2form} as the
map $R$. We will call $R,R'$ {\em braid isomorphisms}. 

\begin{example}\label{ex:doublefusion}
Let us apply internal fusion to the double $D(G)$. 
Fusing the two $G$-factors we get a q-Hamiltonian $G$-space 
$\D(G):=D(G)_{12}$ which is just $G^2$, with $G$-action 
\begin{equation*}
(a,b)^g= (\Ad_g a,\Ad_g b),
\end{equation*}
moment map 
\begin{equation*}
\mu(a,b) = aba^{-1}b^{-1} \equiv [a,b],
\end{equation*}
and 2-form 
\begin{equation*}
\om=\f{1}{2}(a^* \theta,\,b^* \olt) + \f{1}{2} (a^* \bt,\,b^*
\theta)+\f{1}{2}((ab)^*\theta,(a^{-1}b^{-1})^*\olt).
\end{equation*}
The braid isomorphisms in this case are given by 
\begin{equation*}
R:\D(G)\to \D(G),\ (a,b)\mapsto (ab^{-1}a^{-1},\ ab^2)
\end{equation*}

\begin{equation*}
R':\D(G)\to \D(G),\ (a,b)\mapsto (a^{-1}b^{-1}a,\ b^2a) 
\end{equation*}
\end{example}

\begin{remark}
It is interesting to re-examine q-Hamiltonian reduction in connection
with fusion.  Suppose $M$ is a q-Hamiltonian $G\times G\times
H$-space, with moment map $(\mu_1,\mu_2,\mu_3)$. Suppose that $e$ is a
regular value of $\mu:=\mu_1\mu_2$. Then the diagonal action of $G$ on
$\mu^{-1}(e)$ is locally free, hence $M\qu\diag(G):=
\mu^{-1}(e)/\diag(G)$ is an $H$-equivariant orbifold and $\mu_3$
descends to $M\qu\diag(G)$.  We claim that the pull-back of $\om$ to
$\mu^{-1}(e)$ is basic and that the induced form $\om\qu\diag(G)$
makes $M\qu\diag(G)$ into a q-Hamiltonian $H$-space. In fact, since
the pull-back of the additional term in \eqref{Fomega} to
$\mu^{-1}(e)$ vanishes, the claim follows immediately from
\begin{equation*} 
M\qu\diag(G)=(M_{12})_e.
\end{equation*}
Let us also observe that just 
as in the category of Hamiltonian $G$-spaces there is a shifting-trick
for q-Hamiltonian reduction: If $(M,\A,\om,\mu_1,\mu_2)$ is a 
q-Hamiltonian $G\times H$-space, then $f$ is a regular value for 
$\mu_1$ if and only if the identity $e$ is a regular value 
for the moment map on $M\fus \Co^-$ where $\Co=\Ad(G)\cdot f$, 
and in this case 
there is a canonical isomorphism
\begin{equation*}
M_f\cong (M\fus \Co^-)_e=(M\times \Co^-)\qu \diag(G).
\end{equation*}
\end{remark}

\section{Cross-sections and convexity}
One of the basic tools in the study of Hamiltonian $G$-spaces 
is the cross-section Theorem of Guillemin-Sternberg, which
is a method of reducing problems to subgroups of $G$. 
 
In this section we prove a cross-section Theorem for q-Hamiltonian 
$G$-spaces and explain its relation to convexity theorems for the moment 
map.

Let $(M,\A,\om,\mu)$ be a q-Hamiltonian $G$-space and let
$f\in G$. Since the centralizer $Z_f\subset G$ 
is transversal to the conjugacy class $\Co=\Ad(G)\cdot f$ 
there exists 
an open $Z_f$-invariant subset $U\subset Z_f$ containing $f$ and 
and an equivariant diffeomorphism 
$$G\times_{Z_f}U\to G,\,[g,u]\mapsto g\,u$$
onto an open subset of $G$. 
By equivariance of $\mu$, the pre-image $\mu^{-1}(Z_f)$ is a smooth 
$Z_f$-equivariant submanifold $Y\subset M$, and there is a 
natural diffeomorphism $G\times_{Z_f}Y\to M$ onto an open subset. 
In analogy 
to the cross-section theorem of Guillemin-Sternberg
we have
\begin{prop}[Cross-section theorem]\label{propcross}
Let $(M,\A,\om,\mu)$ be a q-Hamiltonian $G$-space, and $f\subset G$, 
$U\subset Z_f$ as above. Then the {\em cross-section} 
$Y:=\mu^{-1}(U)$ is a smooth 
$Z_f$-invariant submanifold, 
and is a q-Hamiltonian $Z_f$-space with the restriction of $\mu$ as 
a moment map. In particular, if $Z_f$ is abelian the cross-section 
is symplectic.
\end{prop}

\begin{proof}
All conditions for a q-Hamiltonian $Z_f$-space are immediate 
except the non-degeneracy condition (B3). 
Let $\iota:\,Y\to M$ be the inclusion. For all
$y\in Y$ the tangent space $T_{\mu(y)}G$ splits into a direct sum 
$$T_{\mu(y)}G=T_y\Co\oplus \{v_\xi(\mu(y))|\,\xi\in\z_f^\perp\}.$$ 
Consequently 
\begin{equation}\label{dirsum}
T_yM=T_yY\oplus\{v_\xi(y)|\,\xi\in\z_f^\perp\}.
\end{equation}
The second summand is mapped under $\d_y\mu$ to a subspace 
of the tangent space to conjugacy class $\Ad(G)\cdot{\mu(y)}$.
If $v\in T_yY$, $\xi\in\z_f^\perp$ we have 
\begin{equation*}
\om(v_\xi,v)=\frac{1}{2} \iota(v)\mu^*(\theta+\olt,\xi)=0
\end{equation*}
because $\iota(d_y\mu(v))(\theta+\olt)\in \z_f$, 
so that the decomposition \eqref{dirsum} is $\om$-orthogonal.
Thus if 
$v\in \ker\iota^*\om$ then also $v\in \ker\om$. Using that $M$ satisfies 
(B3) this shows $v=v_\xi$ 
for some $\xi\in \z_f$ satisfying $\Ad(\mu(x))\xi=-\xi$. 
\end{proof}

\begin{remark}
Instead of the restriction $\mu|Y$, one can also use 
the shifted moment map $\hat{\mu}=f^{-1}(\mu|Y)$. It satisfies 
$\hat{\mu}(x)=e$ for all $x\in\mu^{-1}(f)$.  Remark \ref{rem:qhamham}
applies and shows that $Y$ is locally equivalent near $\mu^{-1}(f)\subset Y$ 
to a Hamiltonian $Z_f$-space in the usual sense. 
\end{remark}

Suppose now that the group $G$ is in addition connected and simply
connected. Then there are canonical choices for the cross-sections
constructed as follows.  Let $T\subset G$ be a maximal torus of $G$,
with Lie algebra $\t$, and $\t_+$ some choice of a positive Weyl
chamber. Every $\Ad(G)$-orbit in $\g$ passes through a unique point of
$\t_+$ so that $\t_+=\g/\Ad(G)$. Let $\Alc\subset\t_+$ be the
fundamental Weyl alcove. Every conjugacy class $\Co\subset G$ contains
a unique point of $\exp(\Alc)\subset T$ so that we can identify
\begin{equation*}
\Alc=G/\Ad(G)
\end{equation*}
as the space of conjugacy classes. For every open face $\sig\subset \Alc$
the centralizer $Z_{\exp(\xi)}$ with $\xi\in\sig$ is independent of 
$\xi$ and will be denoted $Z_\sig$. For $\sig$ in the interior 
of $\Alc$ we have $Z_\sig=T$. Introducing a partial order 
be setting $\sig\succeq\tau$ if $\ol{\sig}\supseteq\tau$ we 
have $\sig\succeq\tau\Rightarrow Z_\sig\subset Z_\tau$, in particular 
every $Z_\sig$ contains $T$. Let us write 
\begin{equation*}
\Alc_\sig=\bigcup_{\tau\succ\sig}\tau
\end{equation*}
and 
\begin{equation*} 
U_\sig=\Ad(Z_\sig)\exp(\Alc_\sig).
\end{equation*}
Then $U_\sig\subset Z_\sig\subset G$ is smooth, and is a slice for the 
$\Ad(G)$-action
at points in $\sig$. In particular, for every $g\in U_\sig$
we have $\Ad(G)g\cap U_\sig=\Ad(Z_\sig)g$. 

\begin{remark}\label{GoldmanFlows}
Let us consider in particular the cross-section $Y=Y_\sig$ for $\sig=\on{int}\Alc$.
The $T$-action over $Y$ extends by equivariance to a $G$-equivariant
$T$-action over the open subset $G\cdot Y\subset M$. 
In fact, it is Hamiltonian in the sense of Proposition 
\ref{HamDyn}. To see this let $q:\,G\to G/\Ad(G)\cong \Alc$ 
be the quotient map. Note that $q$ is smooth over $G\cdot U_\sig$ 
for  $\sig=\on{int}\Alc$. Using \ref{propcross}, it is clear that 
the components of $q\circ \mu$ generate the $T$-action just described.
In the case of moduli space these $T$-actions are known as the 
Goldman flows \cite{G}. See e.g. 
\cite{MW1} for a discussion and references.
\end{remark}
 
The above fact that the cross-sections are equivalent, after shift 
of the moment map, to Hamiltonian spaces in the usual sense implies
in  particular that Kirwan's theorem
on convexity and connectedness of fibers of moment maps applies
to every connected component of the cross-section. 

Using an argument as in \cite{MW1} we can show that the cross-sections
for a connected, simply connected compact Lie group are necessarily 
connected, 
which  then has the following consequence.
\begin{theorem}[Convexity Theorem]\label{th:conv} 
Let $(M,\Alc,\om,\mu)$ be a connected q-Hamiltonian $G$-space,
where $G$ is compact, connected and simply connected. Then 
all fibers of the moment map $\mu$ are connected, and the 
intersection $\mu(M)\cap\Alc$ is a convex polytope.
\end{theorem}
We will not give the detailed argument here since, as the 
following section shows, it is indeed just a translation 
of the argument given in \cite{MW1} into the terminology
of q-Hamiltonian $G$-spaces. Note that the fibers of $\mu$ are 
not necessarily connected if $G$ is for instance a torus.

\section{Relation to Hamiltonian $LG$-spaces}
\label{Loopgroup}
In this Section we prove that there exists a one-to-one correspondence
between q-Hamiltonian $G$-spaces and Hamiltonian $LG$-spaces with
proper moment map. So, one always has a choice either to work with
infinite-dimensional objects ($LG$-spaces) and more conventional
definitions (Hamiltonian spaces) or to use finite-dimensional objects
and the new definitions (q-Hamiltonian $G$-spaces).

\subsection{Loop group $LG$: notations}

Let $G$ be a compact Lie group. We define the loop group $LG$ as a
space of maps
\begin{equation*}
LG= \on{Map}(S^1,G)
\end{equation*}
of a fixed Sobolev class $\lambda>1/2$. Then $LG$ consists 
of continuous maps and the group multiplication is defined pointwise.
Its Lie algebra is the space of maps $L\g=\Om^0(S^1,\g)$ of 
Sobolev class $\lambda$. We define $L\g^*$ as the space of 
1-forms 
\begin{equation*}
L\g^*=\Om^1(S^1,\g) 
\end{equation*}
of Sobolev class $\lambda-1$. The natural pairing of $L\g^*$ 
and $L\g$ given by 
\begin{equation*}
\l A, \xi\r = \oint_{S^1} (A, \xi)
\end{equation*}
makes $L\g^*$ into a subset of the topological dual $(L\g)^*$. 

We view $L\g^*$ as the affine space of connections on the trivial bundle 
$S^1\times G$
and let the loop group $LG$ act by gauge transformations:
\begin{equation*} 
A^g= \Ad_g\, A - g^*\ol{\theta}.
\end{equation*}
%
%
%
Let 
\begin{equation*} 
\Hol:\,L\g^*\to G
\end{equation*}
denote the holonomy map. 
Recall that if 
we identify $S^1\cong \R/\Z$ and let $s\in\R$ denote the 
local coordinate, $\Hol=\Hol_1$ where 
$\Hol_s:\,L\g^*\to G$ is defined as the unique solution 
of the differential equation
\begin{equation*}
\Hol_s(A)^{-1}\,\f{\p}{\p s}\Hol_s(A)=A,\ \ \Hol_0(A)=e.
\end{equation*}
On constant connections $A=\xi\d s$ (for $\xi\in\g$) 
the holonomy map 
restricts to the exponential map: $\Hol_s(\xi\d s)=\exp(s\xi)$.
The map $\Hol_s$ satisfies the equivariance 
condition 
\begin{equation}\label{Transformation}
\Hol_s(A^g)=g(0)\,\Hol_s(A)\,g(s)^{-1},
\end{equation}
in particular the holonomy map $\Hol=\Hol_1$ is equivariant
with respect to the evaluation homomorphism $LG\to G,\,g\mapsto g(1)$
and the adjoint action of $G$ on itself.

Let the {\em based loop group} $\Om G \subset LG$ be 
defined as the kernel of 
the evaluation mapping $LG\to G,\, g\mapsto g(1)$. Then $LG$ is a 
semi-direct product $LG=\Om G\rtimes G$. 
The action of $\Om G$ on $L\g^*$ is free, and the quotient map 
is just the holonomy map. Thus $\Hol:\,L\g^*\to G$ is the universal 
$\Om G$-principal bundle. 
Consider now the closed three-form $\gen\in\Om^3(G)$. Since 
$L\g^*$ is an affine space the pull-back $\Hol^*\gen$ to 
$L\g^*$ is exact. An 
explicit potential for $\Hol^*\gen$ is given as follows.
Define the following 2-form on $L\g^*$: 
\begin{equation*}
\ve=\f{1}{2}\int_0^1 \d s \,(\Hol_s^*\ol{\theta},\,\f{\p}{\p s}\Hol_s^*\ol{\theta}).
\end{equation*}
\begin{prop}\label{veprop}
The 2-form
$\ve$ is $LG$-invariant and its differential is given by 
\begin{equation}\label{veclosed}
\d\ve=-\Hol^* \gen.
\end{equation}
Its contraction with a fundamental vector field $v_\xi$ ($\xi\in L\g$)
is given by 
\begin{equation}\label{vecontr}
\iota(v_\xi)\ve =-\d\oint_{S^1} (A,\xi)+\f{1}{2}\Hol^*
(\theta+\ol{\theta},\xi(0))
\end{equation}
\end{prop}
The proof of this proposition is deferred to the appendix.

\subsection{Equivalence theorem}

\begin{definition}
A Hamiltonian $LG$-space is a Banach manifold $N$
together with  an $LG$-action $\bfA$, 
an invariant 2-form $\sigma\in\Om^2(M)^{LG}$, 
and an equivariant map $\Phi\in C^\infty(N,\Lg^*)^{LG}$ 
such that:
\begin{itemize}
\item[(C1)]
The 2-form $\sigma$ is closed.
\item[(C2)]
The map $\Phi$ is a moment map for the $LG$-action $\bfA$:
\begin{equation*}  
\iota(v_{\xi}) \sigma = d \oint_{S^1} (m,\xi).
\end{equation*}
\item[(C3)]
The form $\sigma$ is weakly non-degenerate, that is, the
induced map $\sigma_x^\flat: T_xN \rightarrow T^*_xN$ is injective.
\end{itemize}
\end{definition}

We will show in this section that every Hamiltonian $LG$-space
with proper moment map determines a q-Hamiltonian $G$-space 
and vice versa. 
Since the action of the based loop group $\OG \subset LG$ on $\Lg^*$ 
is free, its action on $N$ is free as well and we can form the 
quotient 
\begin{equation*}
M=\Hol(N):=N/\OG.
\end{equation*}
If the moment map $\Phi$ is proper then $\Hol(N)$ is a 
smooth finite-dimensional manifold.
We denote by $\Hol$ the projection  $\Hol:\,N\rightarrow \Hol(N)$. 
Since $G=LG/\OG$ 
the diagram

\begin{equation}\label{Eaction}
\vcenter{\xymatrix{
{LG \times N}\ar[d] \ar[r] & N\ar[d]\\
{G\times \Hol(N)}\ar[r] & \Hol(N)
}}
\end{equation}

defines a $G$-action $\Hol(\bfA)$ on $\Hol(N)$ and 
the diagram

\begin{equation}\label{Emoment}
\vcenter{
\xymatrix{
{N}\ar[d] \ar[r] & L\g^*\ar[d]\\
\Hol(N)\ar[r] & G
}}
\end{equation}

a $G$-equivariant map $\Hol(\Phi):\,\Hol(N)\to G$.
The following result says that holonomy manifolds of Hamiltonian 
$LG$-spaces with proper moment maps carry canonically the structure
of q-Hamiltonian $G$-spaces.

\begin{theorem}[Equivalence Theorem]\label{th:equivalence}
Let $(N, \bfA, \sigma, \Phi)$ be a Hamiltonian $LG$-space with 
proper moment map
and let $M=\Hol(N)$ be its holonomy manifold, with 
$G$-action $\A=\Hol(\bfA)$ and map $\mu:=\Hol(\Phi)$. 
The 2-form on $N$
\begin{equation} \label{tO}
\sigma + \Phi^* \ve
\end{equation}
is basic with respect to the projection $\Hol:\,N\rightarrow M$,
and is therefore the pull-back $\Hol^*\om$ of a unique 
2-form $\om$ on $M$. Then $(M, \A, \om, \mu )$ is
a q-Hamiltonian $G$-space. Conversely, given a 
 q-Hamiltonian $G$-space $(M, \A, \om, \mu)$ there is a unique 
Hamiltonian $LG$-space $(N,\bfA,\sigma,\Phi)$ such that $M$ 
is its holonomy manifold.
\end{theorem}

\begin{proof}
Since $\ve$ is $LG$-invariant the 2-form \eqref{tO} is $LG$-invariant.
Moreover if $\xi\in L\g$, we have by \eqref{vecontr}
\begin{eqnarray}\label{tiomcontr}
\iota(v_\xi)(\sigma + \Phi^* \ve)&=&\d\l\Phi,\xi\r+\Phi^*\iota(v_\xi)\ve\\
&=&\f{1}{2}\Phi^*\Hol^*(\theta+\olt,\xi(0))
=\f{1}{2}\Hol^*\mu^*(\theta+\olt,\xi(0)).\nonumber
\end{eqnarray}
This shows that $\sigma + \Phi^* \ve$ is basic, and that the 2-form 
$\om$ on $M$ defined by it satisfies condition (B2). 
Condition (B1) is a consequence of $\d\sigma=0$ and
\eqref{veclosed}:
\begin{equation*}
\Hol^*\d\om=\d\Hol^*\om=\d\sigma+
\Phi^*\d\ve=-\Phi^*\Hol^*\gen=-\Hol^*\mu^*\gen.
\end{equation*}
It remains to check the non-degeneracy condition (B3). 
The kernel of $\om$ is the push-forward by the tangent map $d\Hol$ of 
the kernel of the form $\sigma+\Phi^*\ve$. Suppose $v\in T_yN$ is in the 
kernel. Then 
\begin{equation*}
\iota(v)\sigma_y = -\iota(v)\Phi^* \ve.
\end{equation*}
Since the 1-form on the right hand side annihilates the kernel of $\d_y\Phi$ 
this equation says $v\in \ker(d_y\Phi)^\sigma$. Using cross-sections for 
Hamiltonian $LG$-actions \cite{MW1,MW2} one sees that  
$\ker(d_y\Phi)^\sigma=T_y(LG\cdot y)$, just as for finite dimensional 
Hamiltonian $G$-spaces. Hence there exists $\xi\in L\g$ with 
$v=v_\xi(y)$. By \eqref{tiomcontr} we arrive
at the condition
\begin{equation*}
\Phi^* \Hol^* (\theta+\bt, \xi(0))=0.
\end{equation*}
Applying an arbitrary vector field $v_{\eta}$ to the left hand side
one finds that 
\begin{equation*}
\xi(0)\in \ker(\Ad_{f^2}-1)=\ker(\Ad_{f}-1)\oplus \ker(\Ad_{f}+1)
\end{equation*}
where $f:=\Hol(\Phi(y))$. For every $\xi$ with 
$\xi(0)\in\ker(\Ad_{f}+1)$, the vector $v_{\xi}$ is in the kernel
of $\ti{\sig}$ and the projection by $\d\Hol$ of such vectors 
is the space $\{v_\eta\in T_{\Hol(y)}M,\,\eta\in\ker(\Ad_{f}-1)\}$. 
It remains to show that if $ \xi(0)\in \ker(\Ad_{f}-1)$ 
and $v_\xi\in\ker\ti{\sig}$ 
then $\xi(0)=0$, so that $v_\xi(y)\in T_y(\Om G\cdot y)$. 
Let $\eta\in L\g$ be defined by 
\begin{equation*}
\eta(s):=\Ad\big(\Hol_s(\Phi(y))\big)\xi(0).
\end{equation*}
Since $\eta-\xi\in\Om\g$ the fundamental vector field $v_\eta(y)$
still lies in $\ker\ti{\sig}_y$. On the other hand, using 
Proposition \ref{veprop} and Lemma \ref{Lemma1} in the Appendix 
the value at $\Phi(y)$ of the 
fundamental vector field $v_\eta$ on $L\g^*$ lies in $\ker\ve_{\Phi(y)}$. 
Therefore $v_\eta\in \ker\sig_y$, and by non-degeneracy of 
$\sig_y$ finally $v_\eta(y)=0$. This completes the proof 
that $(M,\A,\om,\mu)$ is a q-Hamiltonian $G$-space.

Suppose conversely that we are given a q-Hamiltonian $G$-space
$(M,\A,\om,\mu)$. Define $N$ as the pull-back of the 
universal $\Om G$-principal bundle $\Hol:L\g^*\to G$ 
by the map $\mu$. In other words, $N$ is given by the fiber product
diagram  \eqref{Emoment}. Define $\Phi$ as  the upper horizontal arrow 
in \eqref{Emoment}. 
There is a unique $LG$-action $\bfA$ on $N$ such that $\Phi$ is equivariant 
and such that the diagram \eqref{Eaction} commutes. Explicitly, it is 
induced from the $LG$-action on $M\times L\g^*$ given by 
\begin{equation*}
g: (x, A) \rightarrow (x^{g(0)}, A^g).
\end{equation*}
The 2-form $\sigma$ is reconstructed from equation \eqref{tO}:
\begin{equation*}
\sigma= \Hol^* \omega - \Phi^* \ve.
\end{equation*}
where $\Hol$ is the left vertical projection in \eqref{Emoment}.
The above argument, read backwards shows that
$N$  is a Hamiltonian $LG$-space with proper moment 
map, and in fact $M=\Hol(N)$. 
\end{proof}

The most 
basic examples of Hamiltonian $LG$-spaces are provided by 
coadjoint $LG$-orbits $\mathcal{O}$
for the affine action on $L\g^*$, equipped with the Kirillov-Kostant-Souriau 
symplectic
structure. All such orbits are preimages ${\mathcal O}= \Hol^{-1}(\Co)$ of 
conjugacy classes in $\Co\subset G$, and conversely 
the holonomy manifold of the orbit ${\mathcal O}$
is just the  conjugacy class $\Co$. By construction, the corresponding
$G$-action is the restriction of the adjoint action of $G$ and the
moment map $\mu$ is the embedding into $G$.
By Theorem \ref{th:equivalence} the conjugacy class $\Co$
inherits from $\mathcal{O}$ a quasi-Hamiltonian structure. By Proposition 
\ref{pr:Conj} such a structure is unique and is described by 
formula \eqref{qKKS}. We have proved the following proposition: 

\begin{prop}\label{qCoadj}
Let ${\mathcal O}= \Hol^{-1}(\Co)$ be an orbit of gauge action
in $L\g^*$. The corresponding holonomy manifold coincides with
the conjugacy class $\Co \subset G$. The induced quasi-Hamiltonian
structure is given by \eqref{qKKS}.
\end{prop}

\section{Moduli spaces of flat connections on 2-manifolds}\label{Modulispaces}
In this section we apply our techniques to describe the symplectic
structure on the moduli spaces of flat connections on Riemann
surface.
Suppose $\Sig$ is an oriented 2-manifold with non-empty boundary
$\p\Sig=\coprod_{j=0}^r V_j$. The space $\A_{flat}(\Sig)\subset 
\Om^1(\Sig,\g)$ of flat $G$-connections on $\Sig\times G$
is invariant under the action of the gauge group $\mathcal{G}(\Sig)$. 
Let $p_j\in V_j$ be given base points on the boundary components and 
let the restricted gauge group $\mathcal{G}_{res}(\Sig)$ consist 
of gauge transformations that are the identity at the $p_j$'s. 
The quotient  $M(\Sig)=\A_{flat}(\Sig)/\mathcal{G}_{res}(\Sig)$ is a 
smooth, finite dimensional manifold -- in fact $M(\Sig)=G^{2(r+k)}$ 
where $k$ is the genus. It carries a residual $G^{r+1}$ 
action, and the holonomies around the $V_j$ descend to a smooth equivariant 
map $\mu:\,M(\Sig)\to G^{r+1}$. Given a collection of conjugacy classes 
$\Co=(\Co_0,\ldots,\Co_r)\subset G^{r+1}$ the quotient 
\begin{equation}\label{eq:reduc}
M(\Sig,\Co)=\mu^{-1}(\Co)/G^{r+1}
\end{equation}
is the moduli space of flat connections with prescribed holonomies,
which  according to Atiyah-Bott \cite{AB} carries a natural
symplectic structure. 
In this section we construct the q-Hamiltonian structure on the moduli
space of flat connections $M(\Sig)$.
The moduli spaces $M(\Sig,\Co)$ with holonomies
in prescribed conjugacy classes $\Co$ are obtained by q-Hamiltonian
reduction from $M(\Sig)$. The main result of this section is Theorem
\ref{th:equivalence1} which identifies $M(\Sig)$ as a fusion product
of a number of copies of the double $D(G)$.  The upshot is that we
arrive at an explicit finite-dimensional description of the symplectic
form on the moduli space. For similar constructions see e.g.
\cite{J1}, \cite{H}, \cite{KS}, and \cite{GHJW}. For {\em complex}
Lie groups there is an alternative finite dimensional construction 
due to Fock-Rosly \cite{FoR}, using ideas from Poisson geometry.  

\subsection{Gauge-theoretic description}
We start by recalling the gauge-theory construction of 
the symplectic 2-form on $M(\Sig,\Co_0,\ldots,\Co_r)$, following 
Atiyah-Bott \cite{AB}. 
Let $\Sig$ be a compact connected oriented 2-manifold 
with boundary components $V_0,\ldots,V_r$ and  
$P\to \Sig$
a principal $G$-bundle. 
For simplicity we assume that $P$ is the trivial bundle 
$\Sig\times G$ although the following discussion goes through 
for non-trivial bundles as well\footnote{Recall that if $G$ is 
simply connected, every $G$-principal bundle over $\Sig$ is trivial.}.
Fix $\lambda>1$ and let $\Om^j(\Sig,\g)$ denote $\g$-valued differential forms 
of Sobolev class $\lambda-j$. Then the gauge $\Gau(\Sig)=\on{Map}(\Sig,G)$ 
is a Banach Lie group modeled on $\on{Lie}(\Gau(\Sig))=\Om^0(\Sig,\g)$.
Consider the space of connections $\A(\Sig)=\Om^1(\Sig,\g)$ as a 
an affine Banach space, with smooth $\Gau(\Sig)$-action given by 
\begin{equation*}
A^g=\Ad(g)A-g^*\olt.
\end{equation*} 
The space $\A(\Sig)$ carries a natural symplectic form 
$$  \sig(a,b)=\int_\Sig a\wedge b ,\ \ \ \
\mbox{$a,b \in T_A\A(\Sig)\cong
\Om^1(\Sig,\g)$}$$
for which the $\Gau(\Sig)$-action is Hamiltonian, with moment map 
\begin{equation}\label{MMoment}
\l\Psi(A),\xi\r=\int_\Sig (\curv(A),\xi)+\int_{\p\Sig}
   (A,\xi),\ \ \ \ \xi\in\Om^0(\Sig,\g)
\end{equation}
here $\curv(A)\in\Om^2(\Sig,\g)$ denotes the curvature
and we have taken the orientation on $\p\Sig$ opposite to the induced 
orientation. Let us choose a base point 
$p_j$
for every boundary component $V_j$,  and let $\Hol^j(A)\in G$ denote the holonomy of $A$ along the 
loop based at $p_j$ and winding once around $B_j$. 
Given conjugacy classes
$\Co_0,\ldots,\Co_r\subset G$ let 
\begin{equation}\label{eq:modulispace}
M(\Sig,\Co_0,\ldots,\Co_r)=\{A\in\A(\Sig)|\curv(A)=0,\
\Hol^j(A)\in \Co_j\}/\Gau(\Sig)
\end{equation}
be the moduli space of flat connections with holonomies in $\Co_j$.
Since the holonomy of a $G$-connection on $S^1$ is determined up to
conjugacy its gauge equivalence class $M(\Sig,\Co_0,\ldots,\Co_r)$ is
a symplectic quotient.  It does not depend on the choice of Sobolev
class $\lambda$ and is finite-dimensional and compact, but sometimes
singular.  (For technical details, see e.g. \cite{AB,DK}.)

Let us suppose for the rest of this section 
that $r\ge 0$, i.e. that $\Sig$ has at least one boundary 
component. This assumption can be made without loss of generality 
because if $\p\Sig=\emptyset$ and $\hat{\Sig}$ is the 2-manifold obtained
from $\Sig$    
by removing a disk, there is a natural identification 
$M(\Sig)\cong M(\hat{\Sig},\{e\})$. 

The condition $\p\Sig\not=\emptyset$ implies that the subset
$\A_{flat}(\Sig)$ of flat connections is a smooth Banach submanifold
of finite codimension. Let 
$$ \iota:\, \A_{flat}(\Sig)\hra \A(\Sig)$$
denote the embedding. Define the restricted
gauge group $\Gau_{res}(\Sig)$ as the kernel of the evaluation mapping
\begin{equation}\label{eq:qhamred}
\Gau({\Sig})\to G^{r+1},\, g\mapsto (g(p_0),\ldots,g(p_r)).
\end{equation}
In other words $\Gau_{res}(\Sig)$ consists of gauge transformations 
that are the identity at the given base points. Set 
\begin{equation*}
M(\Sig):=   \A_{flat}(\Sig)/\Gau_{res}(\Sig).
\end{equation*}
Since the action of $\Gau_{res}(\Sig)$ is free, 
it is not very hard to 
check that $M(\Sig)$ is a finite dimensional, smooth manifold --
in fact we will see below that it is isomorphic to $G^{2(k+r)}$ where 
$k$ is the genus of $\Sig$. From the
original $\Gau(\Sig)$-action we have a residual action $\A$ of 
$\Gau(\Sig)/\Gau_{res}(\Sig)\cong G^{r+1}$, 
and the collection of  holonomy maps $\Hol^j(A)$ 
descends to an equivariant map $\mu:\,M(\Sig)\to G^{r+1}$.
We have 
\begin{equation*}
M(\Sig,\Co_0,\ldots,\Co_r)=\mu^{-1}(\Co_0,\ldots,\Co_r) /G^{r+1}.
\end{equation*}

\begin{theorem}\label{th:Mod1}
There is a natural $G^{r+1}$-invariant 
2-form $\om$ on $M(\Sig)$ for which 
the quadruple $(M(\Sig),\A,\om,\mu)$ is a 
q-Hamiltonian $G^{r+1}$-space.
\end{theorem}

\begin{proof}
Choose orientation and base point 
preserving parametrizations $V_j\cong S^1$ of the boundary circles  and let
\begin{equation*}
R_j:\,\A(\Sig)\to \Om^1(S^1,\g)\cong L\g^*
\end{equation*}
denote the restriction mapping to the $j$th boundary component.
For any $g\in \Gau_{res}(\Sig)$ 
the restriction to the $j$th boundary
component is in the based loop group $\Om G$. Equation \eqref{MMoment} 
and Proposition \ref{veprop} show that pull-back $\iota^*\ti{\sig}$
to $\A_{flat}(\Sig)$ of the 2-form 
\begin{equation} \label{ts}
\ti{\sig}=\sig+\sum_{j=0}^r R_j^*\ve 
\end{equation}
on $\A(\Sig)$ is $\Gau_{res}(\Sig)$-basic, and is therefore the 
pull-back of a unique invariant 2-form $\om\in\Om^2(M(\Sig))$
which satisfies (B1) and (B2). The 
non-degeneracy condition (B3) is verified along the lines 
of the proof of Theorem \ref{th:equivalence}. Let $\pi:\,\A_{flat}(\Sig)
\to M(\Sig)$ denote the projection. The kernel 
of $\om$ at $\pi(A)$ is the image under $\d\pi$ of the kernel 
of $\iota^*\ti{\sig}$. 

Suppose 
$v\in T_A \A_{flat}(\Sig)$ is in the kernel of $\iota^*\ti{\sig}$. 
By definition of $\ti\sig$ this implies that 
$\iota(v)\sig_A$ annihilates $\ker\d_A\Psi\subset T_A \A_{flat}(\Sig)$. 
Consequently 
$v\in \ker(\d_A\Psi)^\sig$. The $\sig$-orthogonal complement 
of $\ker(\d_A\Psi)$ is equal to the tangent space to the orbit 
through $A$. Hence $v=v_\xi$  for some $\xi\in \on{Lie}(\Gau(\Sig))$.  
As in the proof of Theorem \ref{th:equivalence} 
this implies $\xi(p_j)\in \on{ker}(\Ad(\Hol(\Psi(p_j))^2)-1)$
for all $p_j$, and 
we conclude $ \xi(p_j)\in \on{ker}(\Ad(\Hol(\Psi(p_j))+1)$, 
proving the non-degeneracy condition.  
\end{proof}

Our next result 
that $M(\Sig,\Co)$ is a Hamiltonian reduction from $M(\Sig)$.

\begin{theorem}\label{th:Mod2}
The moduli space $M(\Sig,\Co_0,\ldots,\Co_r)$ is a q-Hamiltonian
reduction of $M(\Sig)$ corresponding to the conjugacy
classes $(\Co_0,  \dots ,\Co_r)\subset G^{r+1}$.
\end{theorem}

\begin{proof}
Let $\mathcal{O}_j=\Hol^{-1}(\Co_j)\subset L\g^*$ denote the coadjoint
$LG$-orbits corresponding to the conjugacy classes $\Co_j$, and
$\Om_j$ their Kirillov-Kostant-Souriau symplectic forms. Let
$\hat\iota:\,Z\subset \A_{flat}(\Sig)$ denote the set of flat
connections such that $\Hol^j(A)\in\Co_j$, i.e. $R_j^*A\in
\mathcal{O}_j$.  Let $\hat\pi:\,Z\to M(\Sig,\Co_0,\ldots,\Co_r)$ be
the projection.  By definition of the symplectic form $\sig_{red}$ on
$M(\Sig,\Co_0,\ldots,\Co_r)$,
$$
\hat\pi^*\sig_{red}=
\hat\iota^*\Big(\sig- \sum_{j=0}^r R_j^* \Omega_j\Big)=
\hat\iota^*\Big(\ti{\sig} -\sum_{j=0}^r R_j^* (\Omega_j+\ve)\Big)
.
$$
By Proposition \ref{qCoadj},  $\Omega_j+\ve= \Hol^* \omega_j$,
where $\om_j$ is the unique 2-form which defines a quasi-Hamiltonian
structure to $\Co_j$. Since $\Hol\circ R_j=\Hol^j$, this shows
$$
\pi^*\sig_{red}=\iota^*(\ti{\sig}-(\Hol^j)^*\om_j).
$$
By the Remark in Section \ref{shiftingtrick} this Equation shows that $\sig_{red}$
coincides with the 2-form  obtained by q-Hamiltonian reduction
from the space $M(\Sig)$ equipped with the 2-form $\omega$.
\end{proof}

\begin{remark}
Suppose that $G$ is simply connected. Then the restriction mapping
$\Gau(\Sig)\to \Gau(\p\Sig)\cong LG^{r+1}$ is surjective. Let
$\Gau_0(\Sig)$ be its kernel. Then
$\M(\Sig):=\A_{flat}(\Sig)/\Gau_0(\Sig)$ is a Banach manifold and is a
Hamiltonian $LG^{R+1}$-space with proper moment map (see
e.g. \cite{Do,MW1}).  Its holonomy manifold is the finite-dimensional
$G^{r+1}$-space $M(\Sig)$. See \cite{MW2} for applications of
$M(\Sig)$ in this context. 
\end{remark}

\subsection{Holonomy description}
Our next goal is to make the q-Hamiltonian
structure of $M(\Sig)$ more explicit. We start by introducing 
coordinates on $M(\Sig)$.
Choose a system of smooth oriented paths $U_j$ from $p_j$ to $p_0$ 
($j=1,\ldots,r$) on and loops $A_i,B_i$ ($i=1,\ldots,k$) 
based at $p_0$, such that 
\begin{enumerate}
\item 
The paths $U_j,V_j,A_i,B_i$ meet only at $p_0$. 
\item
Letting $\Sig'$ denote the closed 2-manifold obtained from $\Sig$ by
capping off the boundary components, the fundamental group
$\pi_1(\Sig')$ is the group generated by the $A_i,B_i$, modulo the
relation $\prod_{i=1}^k\,[A_i,B_i]=1$.
\item 
The path $V_0$ is obtained by catenation: 
$$ V_0^{-1}=U_1 V_1 U_1^{-1}\cdots U_r V_r U_r^{-1}
\,[A_1,B_1]\cdots [A_k,B_k].$$
\end{enumerate}
Up to $\mathcal{G}_{res}(\Sig)$-gauge equivalence, every flat connection 
on $\Sig$ is completely determined by the holonomies 
$a_i,b_i,u_j,v_j$ along $A_i,B_i,U_j,V_j$ ($i=1,\ldots,k,j=1,\ldots,r$), 
and conversely every 
collection $a_i,b_i,u_j,v_j$ is realized by some flat connection. 
Consequently we have 
$$ M(\Sig)=G^{2(r+k)} $$
with coordinates $(a_i,b_i,u_j,v_j)$. 
The action of $(g_0,\ldots,g_r)\in G^{r+1}$ 
is given by 
\begin{equation}\label{action1}
a_j\mapsto \Ad_{g_0}a_j,\ b_j\mapsto \Ad_{g_0}b_j,\ 
   u_j\mapsto g_0 u_j g_j^{-1},\ v_j \mapsto \Ad_{g_j}v_j
\end{equation}
and the
components of the moment map $\mu$ are 
\begin{eqnarray}\label{mommaps}
\mu_j(a,b,u,v)&=&v_j^{-1},\ \  (j=1,\ldots r)\\
\mu_0(a,b,u,v)&=&\Ad_{u_1}v_1\cdots \Ad_{u_r}v_r
\,[a_1,b_1]\cdots [a_k,b_k]\nonumber.
\end{eqnarray}
We now construct a q-Hamiltonian structure on $M(\Sig)$ as follows. 
Take $r$ copies of the double $D(G)\cong G^2$, with coordinates 
$(u_j,v_j)$ as in Remark \ref{rem:othercoord}, and $k$ copies of 
its ``internal fusion'' $\D(G)\cong G^2$, with coordinates 
$(a_i,b_i)$ as in Example \ref{ex:doublefusion}. 
Recall that $D(G)$ is a q-Hamiltonian $G\times G$-space 
with moment maps $(v_j^{-1},\Ad_{u_j}v_j)$
while 
$\D(G)$ is a $G$-space with moment map $[a_i,b_i]$. 
Fusing the $D(G)$'s with respect to the second component of the 
$G\times G$-action in each copy together with all $\D(G)$'s
we obtain a q-Hamiltonian $G^{r+1}$-space with action  given by
\eqref{action1} and moment map by \eqref{mommaps}.

\subsection{Equivalence of the gauge theory construction and the holonomy 
construction}
We will now prove that the fusion product from the preceeding 
subsection 
gives indeed the correct 2-form.
\begin{theorem}\label{th:equivalence1}
Let $\Sig$ be a smooth 2-dimensional orientable manifold of
genus $k$ with $r+1$ boundary components.
Then the moduli space $M(\Sig)=\A_{flat}(\Sig)/\mathcal{G}_{res}(\Sig)$ 
is canonically isomorphic to the fusion product 
\begin{equation}\label{fusionprod}
\underbrace{D(G)\fus\cdots D(G)}_{r\ \on{times}}
\fus \underbrace{\D(G)\fus\cdots \D(G)}_{k\ \on{times}}.
\end{equation}
\end{theorem}
\begin{proof}
Let $P$ denote 
the polyhedron obtained by cutting $\Sig$ along the 
paths $U_j,A_i$ and $B_i$. The boundary $\p P$ consists of 
$3r+4k+1$ segments:
\begin{equation} \label{pP}
\partial P= U_1 V_1 U_1^{-1} \dots  U_r V_r U_r^{-1}
A_1B_1A_1^{-1}B_1^{-1} \dots A_kB_kA_k^{-1}B_k^{-1}V_0.
\end{equation}
Since $P$ is contractible, every 
flat connection $A$ on $\Sig$ determines a unique 
function $\psi\in\on{Map}(P,G)$ such that 
\begin{equation}\label{eq:psi} 
A=\psi^*\theta,\ \ \psi(p_0)=e 
\end{equation}
%
Choose an 
orientation-preserving parametrization $\p P\cong [0,1]$ 
such that $p_0=0$, 
and let $\psi_s$ denote the value of $\psi|\p P$  at $s\in\p P$.
As before let  
$\iota:\,\A_{flat}(\Sig)\hra \A(\Sig)$ denote 
the inclusion.
\begin{lemma}\label{lem:boundaryintegral}
The pull-back of symplectic form $\sig$ to the submanifold of flat connections is given by 
the formula 
$$\iota^*\sig=\frac{1}{2} \int_0^1 ds (\psi_s^* \bt,
\frac{\partial}{\partial s} \psi^*_s \bt).$$
\end{lemma}
For a proof of this result see  e.g.  \cite{AM}. 

To proceed we introduce some more
notation. Suppose $\Del=[s_0,s_1]\subset [0,1]$ is 
a subinterval, and $c:=\psi_{s_0}$ and $d=c^{-1}\psi_{s_1}$.
Then 
$$\psi^{s_0}:=c^{-1}\psi$$ 
satisfies Equation \eqref{eq:psi} with initial condition 
$\psi(s_0)=e$ instead of $\psi(p_0)=e$. Set $\psi^{s_0}_s=c^{-1}\psi_s$ 
and define
$$ \ve_{\Del}=\f{1}{2}\int_{s_0}^{s_1}\d s \big((\psi^{s_0}_s)^*\olt
,\,\f{\p}{\p s} (\psi^{s_0}_s)^*\olt\big). $$
Then  
\begin{eqnarray}
\frac{1}{2} \int_{s_0}^{s_1}ds (\psi_s^* \bt,
\frac{\partial}{\partial s} \psi^*_s \bt)&=&
\frac{1}{2} \int_{s_0}^{s_1}ds \big((c\psi^{s_0}_s)^* \bt,
\frac{\partial}{\partial s} (c\psi^{s_0}_s)^* \bt\big)\nonumber \\&=&
\ve_{\Del}+\frac{1}{2}(c^*\theta,d^*\olt).\label{eq:crossterm}
\end{eqnarray}
Lemma \ref{lem:boundaryintegral} asserts that 
$\iota^*\sig=\ve_{\p P}$,
while the 2-form $\ti{\sig}$ on $\A(\Sig)$ that gives to the 
q-Hamiltonian structure on the moduli space $M(\Sig)$  
satisfies
\begin{equation}
\label{tils}
\iota^*\ti\sig=\ve_{\p P}-\sum_{j=0}^r\ve_{V_j^{-1}}.
\end{equation}
Using \eqref{eq:crossterm} we evaluate this equation as follows. 
Combine the segments \eqref{pP} of the boundary $\p P$ to the 
following $r+1+k$ loops on $\Sig$: 
$$
\Del_i=\left\{\begin{array}{c@{\quad:\quad}l} 
V_0&i=0\\
A_i\,B_i\,A_i^{-1}B_i^{-1}&1\le i\le k\\
U_{i-k} V_{i-k} U_{i-k}^{-1}&k+1\le i\le k+r
\end{array}\right.$$
For $i=0,\ldots k+r+1$ let  
$t_i\in[0,1]$ such that $\Del_i=[t_i,t_{i+1}]$, and let 
$c_i:=\psi_{t_i}$ and $d_i=c_i^{-1}c_{i+1}$. Then  
\begin{equation}
\ve_{\p P}=\om_{V_0}+
\sum_{i=1}^{k+r} \Big(\ve_{\Del_i}+\f{1}{2}(c_i^*\theta,d_i^*\olt)\Big).
\end{equation}
Let us first consider the contribution
of the loops $\Delta_{k+j}$ for $j=1,\ldots,r$. By another
application of \eqref{eq:crossterm} we have
\begin{eqnarray*}
\ve_{\Delta_{k+j}} &=& \ve_{U_j} + \ve_{V_j} + \ve_{U_j^{-1}} +
\frac{1}{2}(u_j^* \theta, (v_ju_j^{-1})^* \bt ) +
 \frac{1}{2}(v_j^*\theta, (u_j^{-1})^* \bt ) \} \\
&=& \ve_{V_j} + \frac{1}{2} 
(\Ad_{v_j} u_j^*\theta, u_j^*\theta)+
\frac{1}{2}(u_j^*\theta, v_j^*\theta+v_j^*\bt ) \} .
\end{eqnarray*}
In the first line contributions $\ve_{U_j}$ and  $\ve_{U_j^{-1}}$
cancel each other due to the difference in orientation. 
The term $\ve_{V_j}$ cancels the corresponding contribution
to \eqref{tils}. The remaining expression gives the 2-form on the 
double $D(G)$.

Next, we analyze the contribution of
$\Delta_i=A_i\,B_i\,A_i^{-1}B_i^{-1}$ for $i=1,\ldots,k$:
\begin{eqnarray*}
\lefteqn{\ve_{\Delta_i} = \ve_{A_jB_j} +  \ve_{A_j^{-1}B_j^{-1}}
+ \frac{1}{2}((a_jb_j)^*\theta, (a_j^{-1}b_j^{-1})^*\bt )} \\
&=& \ve_{A_j} + \ve_{B_j} + \ve_{A_j^{-1}} + \ve_{B_j^{-1}}
+ \frac{1}{2} \big( (a_j^*\theta,b_j^*\bt )+ (a_j^*\bt, b_j^*\theta)
+ ((a_jb_j)^*\theta, (a_j^{-1}b_j^{-1})^*\bt ) \big) \\
&=&  \frac{1}{2} \big( (a_j^*\theta,b_j^*\bt )+ (a_j^*\bt, b_j^*\theta)
+ ((a_jb_j)^*\theta, (a_j^{-1}b_j^{-1})^*\bt ) \big) .
\end{eqnarray*}
The last line reproduces the 2-form on ${\D}(G)$.  Finally, the
contribution $\ve_{V_0}$ cancels the corresponding term in
\eqref{tils}.  We have shown that the 2-form $\tilde{\sig}$ is a
pull-back of the sum of $r$ copies of the 2-form on $D(G)$, along the
maps $u_i, v_i$, and of $k$ copies of the 2-form on ${\mathbf D}(G)$,
along the maps $a_i,b_i$. The cross terms $\f{1}{2}(c_i^*\theta,d_i^*\olt)$
are precisely the extra terms coming from fusion \eqref{Fomega}. This
completes the proof of Theorem \ref{th:equivalence1}.
\end{proof}

Suppose now that we are given a tuple of conjugacy classes
$\Co=(\Co_0,\ldots,\Co_r)$. 
Recall that the reduction of $D(G)$ at $\Co_i$ is equal to
$\Co_i^-$. Theorem \ref{th:equivalence1} together with Theorem
\ref{th:Mod2} show that the moduli space $M(\Sig,\Co)$ is a
q-Hamiltonian reduction, 
$$ M(\Sig,\Co)=(\Co_0^-\fus \cdots\fus\Co_r^-\fus\D(G)\fus \cdots\fus\D(G))_e.$$
In particular, the moduli space for the sphere with $d$ holes is the 
reduction of a $d$-fold fusion product of conjugacy classes. This fits 
nicely with the well-known similarities of this space with the symplectic 
reductions of a $d$-fold product of coadjoint $G$-orbits. 

\subsection{Action of the mapping class group}
Let $\Diff(\Sig)$ denote the group of orientation preserving
diffeomorphisms of $\Sig$, with the $C^1$-topology.  
Its Lie algebra is the space of vector
fields that are tangent to the boundary.  The action
\begin{equation*}\Diff(\Sig)\times \A(\Sig)\to\A(\Sig),\  
A^\phi=(\phi^{-1})^*A 
\end{equation*}
preserves the 2-form $\sig$, and is Hamiltonian with moment map
 \begin{equation*}
\l\Phi(A),X\r=-\int_\Sig (\curv(A),\iota(X)A)-
              \int_{\p\Sig}(A,\iota(X)A)
\end{equation*}
for $X\in \on{Lie}(\Diff(\Sig))$.
It also acts on the gauge group $\Gau(\Sig)$ by $g^\phi=(\phi^{-1})^*g$, 
and combines with the gauge group action to an action 
of the semi-direct product $\Gau(\Sig)\rtimes \Diff(\Sig)$. 

Consider the smaller group $\Diff_{res}(\Sig)$ consisting of all 
$\phi\in \Diff(\Sig)$ that preserve the base points $\{p_0,\ldots,p_r\}$
up to permutation. Its action descends to $M(\Sig)$, and combines 
with the action of $G^{r+1}$ to an action of the semi-direct product 
$G^{r+1}\rtimes \Diff_{res}(\Sig)$, where $\Diff_{res}(\Sig)$ acts 
on $G^{r+1}$ by permuting factors. Clearly the $\Diff_{res}(\Sig)$-action 
preserves $\om$ and 
permutes the components of $\mu$.
However, 
the action of the identity component 
$\Diff_{res}^0(\Sig)=\Diff_{res}(\Sig)\cap \Diff^0(\Sig)$ on 
$M(\Sig)$ is trivial because every 
$\phi\in  \Diff^0(\Sig)$ is connected to the identity by 
a smooth path $\phi_t$, and because
the fundamental vector field $v_X$ of $X\in\on{Lie}(\Diff_{res}(\Sig))$
is equal to $v_\xi$ for $\xi=-\iota(X)A\in \on{Lie}(\Gau(\Sig))$. 
All that remains is therefore the action of the mapping class group
$\Gamma(\Sig)=\Diff_{res}(\Sig)/\Diff_{res}^0(\Sig)$, and we obtain an 
action of the semi-direct product  

\begin{equation*}
G^{r+1}\rtimes \Gamma(\Sig).
\end{equation*}
preserving $\om$. The action of the ``pure'' mapping class group,
i.e. the kernel of the homomorphism $\Gamma(\Sig)\to
S(p_0,\ldots,p_r)$ to the permuation group, descends to a
symplectomorphism of the reduced spaces $M(\Sig,\Co_0,\ldots,\Co_r)$.

The action of $\Gamma(\Sig)$ can be described explicitly in terms of
coordinates on $M(\Sig)$.

\begin{example}
 According to Theorem \ref{th:equivalence1} the moduli space $M(\Sig)$
for $\Sig=\Sig_0^2$ the 2-holed sphere is just the double $D(G)$
considered in the previous section. The element $a$ is interpreted as
 parallel transport along a path from $p_1\in V_1$ to $p_2\in V_2$,
 while $ab$ is the holonomy around the boundary component $V_2$. The
 map
\begin{equation*}
S:D(G)\to D(G),\ (a,b)\mapsto (a^{-1},b^{-1})
\end{equation*} 
corresponds to a diffeomorphism exchanging $V_1$ and $V_2$: Indeed it
satisfies $S^*\om=\om$, but switches the $G$-factors so that  
$S^*(\mu_1,\mu_2)=(\mu_2,\mu_1)$ and 
$S((a,b)^{(g_1,g_2)})=(S(a,b))^{(g_2,g_1)}$.
Another interesting action is given by 
\begin{equation*}
Q:D(G)\to D(G),\ (a,b) \mapsto (aba,a^{-1}).
\end{equation*} 
This action by $Q$ is equivariant and preserves both the 2-form and the moment 
map. It corresponds to a diffeomorphism which rotates one of the 
boundary circles by $2\pi$ while leaving the other one fixed.  
Since $D(G)$ acts as the identity under diagonal reduction 
(that is, $(M\fus D(G))_e\cong M$) this explains the existence of
the twist automorphism, Theorem \ref{Dehn}. 
\end{example}
\begin{example}
The fusion operation $M_1\fus M_2$ can be viewed 
as a diagonal $G^2$-reduction 
$$\big(M_1\times M_2\times M(\Sig_0^3)\big)\qu G^2.$$
with respect to two of the three boundary circles of
$\Sig_0^3$. Choosing an element of the mapping class group exchanging
these two boundary circles we obtain a q-Hamiltonian isomorphism
$M(\Sig_0^3)\to M(\Sig_0^3)$ exchanging two components of the moment
map.  It descends to a q-Hamiltonian isomorphism $M_1\fus M_2\to
M_2\fus M_1$.  This is the origin for the braid isomorphisms discussed
in Section \ref{fusion}.

In \cite{MW2}, a fusion operation
was introduced for Hamiltonian $LG$-manifolds with proper
moment maps. Letting $\M(\Sig_0^3)$ be the Hamiltonian $LG^3$-manifold 
associated to $\Sig_0^3$, the fusion product of two Hamiltonian 
$LG$-manifolds $\M_1$ and $\M_2$ is the diagonal $LG^2$-reduction 
$$ \M_1\fus \M_2=(\M_1\times \M_2\times  \M(\Sig_0^3))\qu LG^2;$$
it is a Hamiltonian $LG$-space with proper moment map.  (By contrast,
the moment map for the direct product $\M_1\times \M_2$ with diagonal
$LG$-action is not proper and does not have the correct equivariance
property.) The holonomy manifold of $\M_1\fus \M_2$ is 
$$ \Hol(\M_1\fus \M_2)=\Hol(\M_1)\fus \Hol(\M_2).$$
\end{example}

\section{Relation to Poisson-Lie $G$-spaces}
In this Section we establish a connection between the theory 
of Poisson-Lie groups and the theory of q-Hamiltonian
$G$-spaces. Although these two theories are not equivalent
to each other, the definition of a Poisson-Lie $G$-space
can be rewritten in a form very similar to the definition of a
q-Hamiltonian $G$-space. 
\subsection{Poisson-Lie $G$-spaces}
We begin with a short exposition of the theory of Poisson-Lie
$G$-spaces of J.-H. Lu and A. Weinstein \cite{L,LW}.  As for 
q-Hamiltonian spaces, the target of the moment map is a non-abelian
Lie group.

Throughout this section $G$ denotes a connected and simply connected 
compact Lie group and $T\subset G$ a maximal torus. 
The inner product $(\,,\,)$ on    
$\g$ induces a complex-bilinear form, still denoted $(\,,\,)$ on its 
complexification $\g^\C$. We will regard $\g^\C$ as a real Lie algebra, 
and let $G^\C$ be the corresponding Lie group. Let $\n\subset \g^\C$ 
be the sum of root spaces for the positive roots
and $\a:=\sqrt{-1}\,\t$. Write $A=\exp(\a)$ and $N=\exp(\n)$. 
We have
the Iwasawa decompositions 
$$ \g^\C=\g\oplus \a\oplus \n,\ \ G^\C=G\,A\,N=A\,N\,G. $$
The pairing of the subalgebra 
$\a\oplus \n$ with $\g$ given by the imaginary part 
of $(\,,\,)$ 
$$ \l \zeta, \eta\r  = \im ( \zeta, \eta) \ , \ \zeta\in \a\oplus \n \ , \ \eta \in \alg{g}.$$
is nondegenerate, and identifies $\a\oplus \n\cong \g^*$. Let 
$G^*=AN\subset G^\C$ be the corresponding simply connected subgroup.
We denote the projection to the first factor by $\alpha: G^*\to A$.
According to Drinfeld the isomorphism $\g^\C=\g\oplus\g^*$ means that 
$G$ is a  {\em Poisson-Lie group}, and $G^*$ its dual Poisson-Lie group. 
By the Iwasawa decomposition any element 
$G^{\C}$ can be uniquely written as a product of elements of $G^*$ and $G$:
\begin{equation*}
G^{\C}=G^*G.
\end{equation*}
Left-multiplication of $G$ on $G^\C$ induces an action of $G$ on
$G^*=G^{\C}/G$ which is known as the (left) dressing action (this
terminology is due to Semenov-Tian-Shansky).  
We denote the fundamental vector fields for the dressing action by 
$v_\xi^\sharp$. 

We will use the same notation $\theta,\olt$ for the Maurer-Cartan
forms on $G^\C$ and its subgroups. 
This does not lead to ambiguities since the
Maurer-Cartan form on a subgroup of a group is just the pull-back of
the Maurer Cartan-form on the group.  Sometimes we denote the
Maurer-Cartan forms on $G^*$ by $\theta_{G^*},\olt_{G^*}$ for clarity.  
Let 
$$ \chi_\C=\f{1}{12}(\theta,[\theta,\theta])\in \Om^3(G^\C,\C) $$
\begin{definition}[Lu]
A Poisson-Lie $G$-space is a $G$-manifold $(M,\A)$ together 
with a 2-form $\omega\in\Om^2(M)$, 
and an equivariant map $\mu\in C^\infty(M,G^*)$
such that the following conditions are fulfilled:
\begin{itemize}
\item[(D1)]
The form $\omega$ is closed.
\item[(D2)]
For all $\xi\in\g$,
\begin{equation*} 
\iota(v_\xi^\sharp)\omega 
= 2 \mu^*\,\langle\olt_{G^*} , \xi\rangle .
\end{equation*}
\item[(D3)]
The form $\omega$ is non-degenerate.
\end{itemize}
The map $\mu$ is called a Poisson-Lie moment map.
\end{definition}

The factor of 2 is introduced into (D2) in order to simplify the
comparison of the definition of a Poisson-Lie $G$-space to the
definition of a q-Hamiltonian $G$-space with $P$-valued moment map.
\begin{remark}
\begin{enumerate}
\item
The 2-form $\om$ is not invariant. Rather, the point 
of the definition is that the action map $\A$ becomes a Poisson map
\cite{L}.
\item
Just as for q-Hamiltonian $G$-spaces, there is a ring structure on the 
category of Poisson-Lie spaces. Given two Lie-Poisson 
$G$-spaces $M_1,\,M_2$ there exists the structure of a Lie-Poisson
$G$-space on $M_1\times M_2$ with symplectic form the sum 
$\om_1+\om_2$ and moment map 
the pointwise product $\mu_1\cdot\mu_2$.
The $G$-action is not simply the diagonal $G$-action but is ``twisted''.
See e.g. \cite{FR}.   
\item
The moment map for Poisson-Lie $G$-spaces has the properties 
\begin{equation}\label{eq:plmoment} 
\im(\mu^*\theta)_x=\g_x^0,\ \ \ker(\d_x\mu)^\om=\{v_\xi^\sharp(x),
\,\xi\in\g\}.
\end{equation}
The proof is analogous to that for Hamiltonian $G$-spaces. 
\end{enumerate}
\end{remark}

\subsection{q-Hamiltonian $G$-space with $P$-valued moment maps}
As we explained in Remark \ref{rem:equivcoh}, the possibility 
of choosing $X=\g^*$ or $X=G$ as target space for a generalized moment 
map relies on the existence of a natural equivariantly closed 
equivariant 3-form $\chi_G$ on $X$. 
In this subsection we present another example of a target $X$ 
with this property.

Let $\tau:\,G^\C\to G^\C$ denote the Cartan involution,
defined by exponentiating the complex conjugation mapping 
$ \g^\C\to \g^\C$,  
and let $I:\,G^\C\to G^\C$ denote the map $Ig=\tau(g^{-1})$. 
We will also use the notation $Ig=g^\dagger$ since for $G= \on{SU}(N)$ 
the complexification is 
$G^\C=\on{Sl}(N,\C)$ and the map $I$ is Hermitian conjugation.
Let $P$ denote the symmetric space 
$$ P=\{g\in G^\C,\,g=g^\dagger\} .$$
The adjoint action of $G$ on $G^\C$ leaves $P$ invariant. Let 
$p:\,P\hra G^\C$ denote the embedding. 
Set  
$$\theta_P = p^* \theta,\ \bt_P=p^* \bt\in\Om^1(P,\g^\C)$$
and define a 3-form $\gen_P\in\Om^3(P)$ by 
\begin{equation*}
\gen_P= \frac{1}{12} p^* \im 
(\bt, [\bt, \bt]).
\end{equation*}
For all $\xi\in\g$, the complex conjugate of the 1-form 
$(\theta,\xi)\in\Om^1(G^\C,\C)$ is $-(I^*\olt,\xi)$. Therefore 
the 1-form $(\theta+I^*\olt,\xi)$ is purely imaginary. On 
$P$ the map $I$ is trivial, so that $(\theta_P+\olt_P,\xi)$
is purely imaginary. Put differently, $\theta_P+\olt_P$ takes values 
in $\sqrt{-1}\g$. 
The equivariantly closed extension of $\gen_P$ is the 3-form 
$\gen_{P,G}\in\Om^3_G(P)$ defined by 
$$  \gen_{P,G}(\xi)=\gen_P+\frac{1}{2\sqrt{-1}}\, \mu_P^* 
(\theta_P+ \bt_P, \xi). $$

\begin{definition}
A q-Hamiltonian $G$-space with $P$-valued moment map is a manifold $M$ equipped with a $G$-action
$\A$, a 2-form $\omega_P\in\Om^2(M)$ 
and an equivariant map $\mu_P\in C^\infty(M,P)^G$ such that:
\begin{itemize}
\item[(E1)]
The differential  of $\omega_P$ is given by:
\begin{equation*} 
\d\omega_P = - \mu^*_P \gen_P.
\end{equation*}
\item[(E2)]
For all $\xi\in\g$,
\begin{equation*} 
\iota(v_{\xi}) \omega_P= \frac{1}{2\sqrt{-1}}\, \mu_P^* 
(\theta_P+ \bt_P, \xi).
\end{equation*}
\item[(E3)]
The form $\omega_P$ is non-degenerate.
\end{itemize}
The map $\mu_P$ is called a $P$-valued moment map.
\end{definition}

Since  $G$ by assumptionis connected, equivariance of the moment map
$\mu_P$ together with conditions (E1) and (E2) imply invariance of the
2-form $\om_P$:
\begin{eqnarray*} 
\mathcal{L}_{v_\xi}\om_P&=&\iota(v_\xi)\d\om_P+\d\iota(v_\xi)\om_P\\
&=&\mu_P^*\Big(-\iota(v_\xi)\chi_P+\frac{1}{2\sqrt{-1}}\,  
\d(\theta_P+ \bt_P, \xi)\Big)=0.
\end{eqnarray*} 

\subsection{Equivalence of Poisson-Lie G-spaces 
and q-Hamiltonian $G$-spaces with $P$-valued moment map} In this
subsection we prove the equivalence of the definitions of a
Poisson-Lie $G$-space and of a q-Hamiltonian $G$-space with $P$-valued
moment map. Consider the map
\begin{equation*}
j: G^\C \rightarrow P \ , \ j(b)=b\, b^\dagger .
\end{equation*}
It turns the left $G$-action on $G^\C$ into the adjoint action on
$P$ and restricts to an equivariant diffeomorphism $G^*\cong P$.
Let $\kappa:\,G^*\hra G^\C$ denote the embedding.

\begin{prop}\label{PLqHamP}
Let $(M, \A, \omega, \mu)$ be a Poisson-Lie $G$-space.
Then the manifold $M$ equipped with the same $G$-action
$\A$, moment  map 
\begin{equation*}
\mu_P = j\circ \mu : M\rightarrow P
\end{equation*}
and 2-form
\begin{equation*}
\omega_P= \omega + \frac{1}{2}
\mu^*\kappa^* \im \ (I^*\olt, \theta)
\end{equation*}
is a q-Hamiltonian $G$-space with $P$-valued moment map.
\end{prop}

\begin{proof}
The map $\mu_P$ is equivariant because it is the composition of two
equivariant maps. 
To check (E1) observe first that 
\begin{equation}\label{eq:1} 
j^*\gen_\C=\gen_\C+I^*\gen_\C-\f{1}{2}\d(I^*\olt,\theta).
\end{equation}
Taking imaginary parts this shows $j^*\im\gen_\C=-\f{1}{2}
\im\d(I^*\olt,\theta)$. Therefore  
\begin{equation*}
\d \omega_P = \d \omega + \frac{1}{2} \d \mu^* \im (I^*\olt, \theta)
= - \mu^* j^* \im \gen_\C
=- \mu_P^* \gen_P
\end{equation*}
in accordance with (E1). 
By Lemma \ref{lem:j} below the imaginary part of the 
1-form $\iota(v_\xi^\sharp)\kappa^*(\theta,I^*\olt)$ on $G^*$ is given by
$$ \iota(v_\xi^\sharp)\kappa^*\im(\theta,I^*\olt)
=4\kappa^*\im(\olt,\xi)-\kappa^*j^*\im (\theta+\olt,\xi).$$
Using this fact together with (D2) we compute  
\begin{eqnarray*}
\iota(v_\xi)( \omega - \frac{1}{2}
\mu^* \im \ \kappa^*(I^*\olt, \theta))&=&
\f{1}{2} \mu^*\kappa^*j^*\,\im \big(
(\theta+\olt,\xi)
\big)
\\&=&\f{1}{2}\mu_P^*\im(\theta+\olt,\xi)
\end{eqnarray*}
which gives (E2). 
To verify (E3) let 
$v\in \ker(\om_P)_x$, that is, 
\begin{equation*}
\iota(v) \omega = -\iota(v)\frac{1}{2}\im
\mu^*\kappa^*\,(I^*\olt_{G^*}, \theta_{G^*}).
\end{equation*}
Since the right hand side of this equation annihilates the kernel 
of $d_x\mu$, this shows that $v\in \ker(d_x\mu)^\om$. By 
\eqref{eq:plmoment} this implies that $v=v_\xi(x)$ for some $\xi\in\g$.
Using (E2) we arrive at the condition 
\begin{equation} \label{tP}
0=\iota(v_{\xi}) \omega_P = \frac{1}{2\sqrt{-1}}\,
\mu_P^* (\theta_P + \bt_P, \xi).
\end{equation}
Contracting this Equation with $v_\eta$ for $\eta\in\g$ 
shows 
\begin{equation*}
(\eta, \Ad_{\mu_P(x)} \xi - \Ad_{\mu_P(x)^{-1}} \xi) =0
\end{equation*}
for all $\eta$, 
or 
$$\xi\in \ker(\Ad_{\mu_P(x)^2}-1)=
\ker(\Ad_{\mu_P(x)}-1)\oplus \ker(\Ad_{\mu_P(x)}+1).$$
Since the eigenvalues of $\Ad_{\mu_P(x)}$ are nonegative real numbers, 
the space 
$\ker(\Ad_{\mu_P(x)}+1)$ is trivial. 
Hence
$\Ad_{\mu_P(x)} \xi =\xi$ . 
For such $\xi$ the above equation becomes  
\begin{equation*}
0= \mu_P^* (\theta_P + \bt_P, \xi)=2\mu_P^*(\theta_P,\xi).
\end{equation*}
By (D2) this equation shows $\iota(v_\xi)\om=0$, and finally 
$v_\xi=0$ by non-degeneracy of $\om$, proving (E3).
\end{proof}

One can easily reverse the argument and show that the
structure of a q-Hamiltonian $G$-space with $P$-valued moment map defines the structure
of a Poisson-Lie $G$-space on the same manifold.
In the proof we used the following Lemma: 

\begin{lemma}\label{lem:j}
The contraction of the fundamental vector field $v_\xi^\sharp$ for the 
dressing action with the 2-form $\kappa^*(I^*\olt,\theta)$ on $G^*$ is 
given by the formula
\begin{equation} \label{eq:4}
\iota(v_\xi^\sharp)\kappa^*(I^*\olt,\theta)=2\kappa^*
(\olt+
I^*\theta,\xi)-\kappa^*j^*(\theta+\olt,\xi).\end{equation}
\end{lemma}
\begin{proof}
We compute: 
\begin{eqnarray}\label{eq:kk} 
\iota(v_\xi^\sharp)\kappa^*(I^*\olt,\theta)&=&
(\iota(v_\xi^\sharp)\kappa^*I^*\olt,\kappa^*\theta)-\nonumber 
(\kappa^*I^*\olt,\iota(v_\xi^\sharp)\kappa^*\theta)\\
&=&(\iota(v_\xi^\sharp)\kappa^*(\theta+I^*\olt),\kappa^*\theta)
-(\kappa^*I^*\olt,\iota(v_\xi^\sharp)\kappa^*(I^*\olt+\theta))\\
&& -(\iota(v_\xi^\sharp)\kappa^*\theta,\kappa^*\theta)
+(\kappa^*I^*\olt,\iota(v_\xi^\sharp)\kappa^*I^*\olt)\nonumber 
\end{eqnarray}
The last two terms in this expression cancel, for the following reason.
It is easy to see that one of them is a complex conjugate of the other.
We will show that the first one is real which ensures the cancellation.
Indeed, the subalgebra $\n$
equals the kernel of the bilinear form $(,)$ 
restricted to $\a\oplus \n$. Hence
$$ (\iota(v_\xi^\sharp)\kappa^*\theta,\kappa^*\theta)=
(\iota(v_\xi^\sharp)\alpha^*\kappa^*\theta,\alpha^*\kappa^*\theta)$$
which is real since the restriction
of $(,)$ to $\a$ is real-valued.
To compute the first two terms in \eqref{eq:kk} observe that 
$$ I^*\olt+\theta=\Ad_{b^\dagger}j^*\olt $$
which gives 
\begin{equation} \label{eq:2}
\iota(v_\xi^\sharp)\kappa^*(\theta+I^*\olt)=
\Ad_{b^{-1}}\xi-\Ad_{b^\dagger}\xi.\end{equation}
Therefore 
$$ \iota(v_\xi^\sharp)\kappa^*(I^*\olt,\theta)
=\kappa^*(\olt+I^*\theta-\Ad_{(b^\dagger)^{-1}}\theta-\Ad_b I^*\olt,\xi).$$
Combining this with 
\begin{equation} \label{eq:3}
j^*(\theta+\olt)=\Ad_{(b^\dagger)^{-1}}\theta+\Ad_bI^*\olt+\olt+I^*\theta
\end{equation}
completes the proof of the Lemma.
\end{proof}

\subsection{Equivalence of q-Hamiltonian $G$-spaces with $P$-valued moment map 
and Hamiltonian $G$-spaces}

Let us note that the space $P$ has two important properties: 
\begin{enumerate}
\item
$P$ is contractible. Hence, the closed 3-form $\gen_P$ is exact. 
\item
The restriction of the exponential map $\exp:\,\g^\C\to G^\C$
to $\sqrt{-1}\,\alg{g}\subset \alg{g}^{\C}$ is invertible
and has image $P$. Let $\varkappa:\,P\rightarrow \sqrt{-1}\,\alg{g}$
be the inverse map, and identify $\sqrt{-1}\,\alg{g}\cong\g^*$ 
by means of the pairing $\im(\,,\,)$. 
\end{enumerate}

Using these facts one can convert a q-Hamiltonian
$G$-space with $P$-valued moment map into a usual Hamiltonian $G$-space. 

\begin{prop}\label{qHamPHam}
\begin{enumerate}
\item
There exists a canonical $G$-invariant 2-form $\tau$ on $P$ such that:
\begin{equation*}
d\tau = \gen_P ,\ \ 
\iota(v_{\xi})\tau = 
d(\varkappa,\xi) - \frac{1}{2\sqrt{-1}}\, (\theta_P + \bt_P,\xi).
\end{equation*}
\item
Let $(M, \A, \mu_P, \omega_P)$ be a q-Hamiltonian $G$-space with 
$P$-valued moment map.
Then the manifold $M$ with the same $G$-action
$\A$, moment map 
\begin{equation*}
\mu = \varkappa \circ \mu_P : M\rightarrow \alg{g}^*
\end{equation*}
and  2-form
\begin{equation*}
\omega= \omega_P + \mu_P^* \tau
\end{equation*}
is a Hamiltonian $G$-space. 
\end{enumerate}
\end{prop}
This was proved in \cite{A}. In fact, the proof is essentially
as that of Proposition \ref{HamqHam} 
since the restriction of the exponential map to $\sqrt{-1}\g$ is 
invertible. 
Propositions \ref{PLqHamP} and \ref{qHamPHam}
reduce the theory of Poisson-Lie $G$-spaces to the usual theory of 
Hamiltonian G-spaces.

\begin{appendix}
\section{Properties of the form $\ve$} 
In this Appendix we give the proof of Proposition \ref{veprop}
concerning the properties of the 2-form 
\begin{equation*}
\ve=\f{1}{2}\int_0^1 \d s \,(\Hol_s^*\ol{\theta},\,\f{\p}{\p
s}\Hol_s^*\ol{\theta})
\end{equation*}
on $L\g^*$.  We will need the following properties of the 
pull-back of $\ol{\theta}$ 
under $\Hol_s$: 
\begin{lemma}\label{Lemma1} (Properties of $\Hol_s^*\olt$.)
Let $\A_g:\,L\g^*\to L\g^*$ the action defined by $g\in LG$.
Then 
\begin{equation}\label{One}
\A_g^*\Hol_s^*\ol{\theta}=\Ad({g(0)})\Hol_s^*\ol{\theta}.
\end{equation} 
The contractions with fundamental vector fields $v_\xi$ (for $\xi\in L\g$) 
are 
\begin{equation}\label{Two}
\iota(v_\xi)\Hol_s^*\ol{\theta} =
\xi(0)-\Ad({\Hol_s(A)})\xi(s).
\end{equation}
For any $\zeta\in L\g^*$, viewed as a constant vector field on $L\g^*$, 
one has  
\begin{equation}\label{Three}
\iota(\zeta)\Hol_s^*\olt=\int_0^s \Ad({\Hol_s(A)})\zeta(u)\, \d u.
\end{equation} 
\end{lemma}

\begin{proof}
For  $h\in G$ let $R_h,\,L_h:\,G\to G$ be the left-resp. right 
multiplication by $h$.  
By the equivariance property (\ref{Transformation}) and right-invariance 
of $\ol{\theta}$, we have
\begin{equation*} 
\A_g^*\Hol_s^*\ol{\theta}=\Hol_s^* L_{g(0)}^* R_{g(s)^{-1}}^*\ol{\theta}
=\Hol_s^* \Ad({g(0)})\ol{\theta}= \Ad({g(0)})\Hol_s^*\ol{\theta}.
\end{equation*} 
By another application of (\ref{Transformation}), 
the tangent map $d_A\Hol_s$ satisfies
\begin{equation*}
(d_A \Hol_s)(v_\xi)= \xi(0)\Hol_s(A)-\Hol_s(A)\xi(s) 
\end{equation*}
which implies (\ref{Two}). For the third identity write 
$\Hol_s(A+t\zeta)=\phi_s(t\zeta)\Hol_s(A)$
so that 
\begin{equation*}
\iota(\zeta)\Hol_s^*\olt=
(d_A \Hol_s)(\zeta) \,\Hol_s(A)^{-1}=\f{\p}{\p t}\Big|_{t=0}
\phi_s(t\zeta).
\end{equation*}
Differentiating the defining identity
\begin{equation*}
\Hol_s(A)^{-1}\phi_s(t\zeta)^{-1}\f{\p}{\p s} \Big(\phi_s(t\zeta)\Hol_s(A)\Big)=
A+t\zeta
\end{equation*}
with respect to $t$ at $t=0$ gives 
\begin{equation*}
\Ad(\Hol_s^{-1})\f{\p}{\p s}\big(\iota(\eta)\Hol_s^*\olt\big)=\zeta
\end{equation*}
Applying $\Ad(\Hol_s)$ and integrating from $0$ to $s$ gives 
\eqref{Three}.
\end{proof}
We now give the proof of Proposition \ref{veprop}.
\begin{proof} 
$LG$-invariance of $\ve$ follows immediately from (\ref{One}).
The differential is computed as follows:
\begin{eqnarray*}
\d\ve&=& \f{1}{4} \int_0^1 \d s\,(\Hol_s^* [\olt,\olt],\,\f{\p}{\p s}\Hol_s^*\olt)
-\f{1}{4} \int_0^1 \d s\,(\Hol_s^*\olt,\,\f{\p}{\p s}\Hol_s^* [\olt,\olt])\\
&=& \f{1}{2}\int_0^1 \d s\,(\Hol_s^* [\olt,\olt],\,\f{\p}{\p s}\Hol_s^*\olt)
-\f{1}{4} \Hol^*([\olt,\olt],\olt)\\
&=& \f{1}{6}\int_0^1 \d s\, \f{\p}{\p s} \Hol_s^*([\olt,\olt],\olt)
-\f{1}{4} \Hol^*([\olt,\olt],\olt) \\
&=& -\f{1}{12}   \Hol^* ([\olt,\olt],\olt)=-\Hol^*\gen.
\end{eqnarray*}
Given $\xi\in L\g$  we have, by partial integration 
\begin{eqnarray*}
\iota(v_\xi)\ve &=&\f{1}{2}\int_0^1 \d s \,(\iota(v_\xi)\Hol_s^* \olt,\,\f{\p}{\p s} 
\Hol_s^* \olt)-\f{1}{2}\int_0^1 \d s \,(\Hol_s^* \olt,\,\f{\p}{\p s} \iota(v_\xi)
\Hol_s^* \olt)\\
&=&\int_0^1 \d s \,(\iota(v_\xi)\Hol_s^* \olt,\,\f{\p}{\p s} 
\Hol_s^* \olt) -\f{1}{2}\Hol^*  ( \olt,\,\iota(v_\xi)\olt)\\
&=& \int_0^1 \d s \,(\xi(0)-\Ad(\Hol_s)\xi(s),\,\f{\p}{\p s} 
\Hol_s^* \olt) -\f{1}{2}\Hol^* (  \olt-\theta,\,  \xi(0))\\
&=&\f{1}{2} (\Hol^* (\theta+ \olt),\xi(0))- \int_0^1 (\xi(s),
\Ad(\Hol_s^{-1})\,\f{\p}{\p s} 
\Hol_s^* \olt).
\end{eqnarray*}
By (\ref{Three}) we have for all $\zeta\in T_AL\g^*\cong L\g^*$
\begin{equation*}
\iota(\zeta) \int_0^1 (\xi(s),\,\Ad(\Hol_s^{-1})\,\f{\p}{\p s} 
\Hol_s^* \olt)=\int_0^1 (\xi(s),\zeta(s))
=\iota(\zeta)\d \oint (A,\,\xi)
\end{equation*}
which concludes the proof.
\end{proof}

\end{appendix}

\end{document}